\documentclass[11pt]{article} \makeatother
\usepackage{amsfonts, amsmath, amssymb}
\usepackage{graphicx}
\usepackage{subfigure}
\allowdisplaybreaks[4]
\newtheorem{theo}{Theorem}[section]
\newtheorem{lemma}[theo]{Lemma}

\newtheorem{definition}[theo]{Definition}

\newcommand{\be}{\begin{equation}}
\newcommand{\ee}{\end{equation}}
\newcommand\bes{\begin{eqnarray}} \newcommand\ees{\end{eqnarray}}
\newcommand{\bess}{\begin{eqnarray*}}
\newcommand{\eess}{\end{eqnarray*}}

\numberwithin{equation}{section}

\begin{document}
\date{}
\setlength{\baselineskip}{17pt}{\setlength\arraycolsep{2pt}

\font\sevenrm=cmr7
\font\sevenit=cmti7
\font\sevenbf=cmbx7
\makeatletter
\def\ps@plain{
\def\@oddhead{\ifnum\thepage=1\hss\baselineskip8pt
\vtop to 0 pt{\vskip-1truecm\hbox
}\vss}\else\hss\fi}\def\@oddfoot{\rm\hfil EJQTDE, 2010 No. 15, p. \thepage}
\def\@evenhead{}\let\@evenfoot\@oddfoot}
\makeatother
\pagestyle{plain}
\thispagestyle{plain}

\begin{center} {\bf\Large Existence, general decay and blow-up of solutions for a viscoelastic  Kirchhoff equation with Balakrishnan-Taylor damping and dynamic boundary conditions  }\\[2mm]

 {\large Gang Li, Biqing Zhu, Danhua Wang  \\[1mm]
{College of Mathematics and Statistics, Nanjing University of
Information Science and Technology, Nanjing 210044, China.\\ E-mail: {\bf brucechu@163.com}.}}\\[1mm]
\end{center}

\begin{abstract}
Our aim in this article is to study a nonlinear viscoelastic Kirchhoff equation with strong damping, Balakrishnan-Taylor damping, nonlinear source and dynamical boundary condition.
Firstly, we prove the local existence of solutions by using the Faedo-Galerkin approximation method combined with a contraction mapping theorem. We then prove that if
the initial data enter into the stable set, the solution globally exists, and if the initial data enter into the unstable set, the solution blows up in a finite time. Moreover, we obtain a general decay result of the energy, from which the usual exponential and polynomial decay rates are only special cases.

\end{abstract}



\section{Introduction }
\setcounter{equation}{0}

In this paper, we are concerned with the following problem:
\begin{equation}\label{1.1}
\left\{ {{\begin{array}{*{20}l} \displaystyle u_{tt}-M(t)\Delta u +\int_{0}^{t}g(t-s)\Delta u(x,s){\rm d}s + \alpha u_{t}-\beta\Delta u_{t} =|u|^{p-2}u , \ & x\in \Omega,  t > 0, \medskip\\
\displaystyle u(x,t) =0, \ & x\in \Gamma_{0},  t > 0, \medskip\\
\displaystyle u_{tt}(x, t) = \displaystyle-M(t)\frac {\partial u}{\partial \nu}(x,t)+\int_{0}^{t}g(t-s)\frac{\partial u}{\partial \nu}(x,s){\rm d}s-\beta\frac {\partial u_{t}}{\partial \nu}(x,t)\medskip\\\quad\quad\quad\quad\ \ -\gamma |u_t|^{m-2}u_{t}(x,t), \ & x \in \Gamma_{1}, t > 0, \medskip\\
u(x, 0) = u_0(x), u_t(x, 0) = u_1(x), \ & x \in  \Omega,
\end{array} }} \right.
\end{equation}
where $M(t)=a+b\|\nabla u(t)\|_2^2+\displaystyle\sigma\int_\Omega\nabla u(t)\nabla u_t(t){\rm d}x$, $\Omega$ is a regular and bounded domain of $\mathbb{R}^{N}$, $(N\geq1) $, $\partial \Omega= \Gamma_{0}\cup\Gamma_1$, $mes(\Gamma_0)>0$, ${\Gamma_0\cap\Gamma_1= \emptyset}$, $\frac {\partial}{\partial \nu}$ denotes the unit outer normal derivative, $m\geq2$, $p>2$ and $a$, $b$, $\sigma$, $\alpha$, $\beta$, $\gamma$ are positive constants. The relaxation function $g$ is a positive and uniformly decaying function. The function $u_0$, $u_1$ are given initial data.

From the mathematical point of view, these problems like \eqref{1.1} take into account acceleration terms on the boundary. Such type of boundary conditions are usually called dynamic boundary conditions. They come from several physical applications and play an important role in different dimension space (see  \cite{AKS1996}, \cite{BEA1976}, \cite{BST1964}, \cite{GOL2006} for more details).

From the theoretical point of view, such as in the one-dimensional case, many existence, uniqueness and decay results have been established (see \cite{DAL1994}, \cite{DAL1996}, \cite{DAL1999}, \cite{DL2002}, \cite{PEL2008}, \cite{PM2004}, \cite{ZH2007} for more details). For example, in \cite{DAL1999} the author considered the following problem
\begin{equation}\label{1.4}
\left\{ {{\begin{array}{*{20}l} \displaystyle u_{tt}-u_{xx}-u_{txx}=0 , \ & x\in (0,L),  t > 0, \medskip\\
\displaystyle u(0,t) =0, \ & t > 0, \medskip\\
u_{tt}(L, t) = -[u_x+u_{tx}](L,t), \ &  t > 0, \medskip\\
u(x, 0) = u_0(x), u_t(x, 0) = v_0(x), \ & x \in  (0,L), \medskip\\
u(L, 0) = \eta, u_t(L, 0) = \mu \ &  t > 0.
\end{array} }} \right.
\end{equation}
By using the theory of B-evolutions and theory of fractional powers, the author proved that problem \eqref{1.4} gives rise to an analytic semigroup in an appropriate functional space and obtained the existence and the uniqueness of solutions. For a problem related to \eqref{1.4}, an exponential decay result was obtained in \cite{DAL1996}, which describes the weakly damped vibrations of an extensible beam. Later, Zang and Hu \cite{ZH2007}, considered the problem
\begin{equation}
\left\{ {{\begin{array}{*{20}l} \displaystyle u_{tt}-p(u_{x})_{xt}-q(u_{x})_x=0 , \ & x\in (0,l),  t > 0, \medskip\\
\displaystyle u(0,t) =0, \ & t \geq 0, \medskip\\
(p (u_{x})_t+q(u_x)(l,t)+k u_{tt}(l,t)) = 0, \ &  t \geq 0, \medskip\\
u(x, 0) = u_0(x), u_t(x, 0) = u_1(x), \ & x \in  (0,l).\nonumber
\end{array} }} \right.
\end{equation}
By using the Nakao inequality, and under appropriate conditions on $p$ and $q$, they established both exponential and polynomial decay rates for the energy depending on the form of the terms $p$ and $q$. In \cite{PM2004} the authors considered the following linear wave equation with strong damping and dynamical boundary conditions
\begin{equation}\label{1.5}
\left\{ {{\begin{array}{*{20}l} \displaystyle u_{tt}-u_{xx}-\alpha u_{txx}=0 , \ & x\in (0,l),  t > 0, \medskip\\
\displaystyle u(0,t) =0, \ & t > 0, \medskip\\
u_{tt}(l, t) = -\varepsilon[u_x+\alpha u_{tx}+ r u_t](l,t), \ &  t > 0,
\end{array} }} \right.
\end{equation}
as an alternative model for the classical ODE, namely
\begin{equation}\label{1.6}
m_1 u''(t)+d_1u'(t)+k_1(t)=0.
\end{equation}
Based on the semigroup theory, spectral perturbation analysis and dominant eigenvalues, they compared analytically these two approaches to the same physical system. Later, Pellicer in \cite{PEL2008} considered \eqref{1.5} with a control acceleration $\varepsilon f \left(u(l,t),u_t(l,t)/\sqrt\varepsilon\right)$ as a model for a controlled spring-mass-damper system and established some results concerning its large time behavior. By using the invariant manifold theory, the author proved that the infinite dimensional system admits a two dimensional invariant manifold where  the equation is well represented by a classical nonlinear oscillations ODE \eqref{1.6}, which can be exhibited explicitly.

In the multi-dimensional cases, we can cite (\cite{GS2006}, \cite{GT1994}, \cite{VIT1999}) for problems with the Dirichlet boundary conditions and (\cite{TOD1998}, \cite{TOD1999}, \cite{TV2005}) for the Cauchy problems. Recently, Gerbi and Said-Houair in (\cite{GH2008}, \cite{GH2011}) studied problem \eqref{1.1} with $M\equiv1$ and without the relaxation function $g$, linear damping $\alpha u_t$. They showed  in \cite{GH2008} that if the initial data are large enough then the energy and the $L^p$ norm of the solution of the problem is unbounded, grows up exponentially as time goes to infinity. Later in \cite{GH2011}, they established the global existence and asymptotic stability of solutions starting in a stable set by combining the potential well method and the energy method. A blow-up result for the case $m=2$ with initial data  in the unstable set was also obtained. For the other works, we refer the readers to  (\cite{BEA1976}, \cite{GOL2006}, \cite{GH2012}, \cite{SW2005}, \cite{WU2011}, \cite{RU2010})  and the references therein.

Many authors have investigated the existence, decay and blow-up results for nonlinear viscoelastic wave equation with boundary dissipation (see \cite{LEV1998}, \cite{LIU1998}, \cite{MES2003}, \cite{MES2006}). For example,
Cavalcanti et al \cite{CAV2002} considered the following equation
\begin{align}
u_{tt}-\Delta u+\int_0^{t}g(t-\tau)\Delta u(\tau){\rm d}\tau+ a(x)u_t+|u|^{\gamma}u=0,\quad (x,t)\in\Omega\times(0,\infty)\nonumber
\end{align}
where $\Omega$ is a bounded domain of $\mathbb{R}^{n}(n\geq1)$ with smooth boundary $\partial \Omega$, $r>0$ and $a: \Omega\rightarrow \mathbb{R_+}$ is a bounded function, which may vanish on a part of the domain. $g$ is a positive nonincreasing function defined on $\mathbb{R_+}$. Under the condition that $a(x)\geq a_{0}>0$
on $\omega\subset\Omega$, with $\omega$ satisfying some geometry restrictions, when $-\xi_{1}g(t)\leq g'(t)\leq-\xi_{2}g(t)$, $t\geq0$, for some positive constants $\xi_{1}$ and $\xi_{2}$, they proved an exponential rate of decay. Later, S. Berrimi and S. A. Messaoudi \cite{BM2004} improved Cavalcanti's result by introducing a different functional, which allowed to weaken the conditions on both $a$ and $g$.
Recently, in \cite{GH2013}, Gerbi and Said-Houair studied problem \eqref{1.1} with $M\equiv1$, without the linear damping $\alpha u_t$ and got the existence, exponential growth results. Jeong et al \cite{JEO2015} concerded the following problem and proved the general energy decay for nonlinear viscoelastic wave equation with boundary damping.
\begin{equation}
\left\{ {{\begin{array}{*{20}l} \displaystyle u''-\Delta u -\alpha \Delta u'+\int_{0}^{t}h(t-\tau)div[a(x)\nabla u(\tau)]{\rm d}\tau =0 , \ & {\rm in}\ \  \Omega\times(0,\infty), \medskip\\
\displaystyle u =0, \ & {\rm on}\ \ \Gamma_{0}\times(0,\infty), \medskip\\
\displaystyle u''+\frac {\partial u}{\partial \nu}+\alpha\frac {\partial u'}{\partial \nu}-\int_{0}^{t}h(t-\tau)[a(x)\nabla u(\tau)]\cdot\nu{\rm d}\tau + |u'|^{m-2}u'=|u|^{p-2}u,\medskip\\\quad\quad\quad\quad\ \  \ & {\rm on}\ \ \Gamma_{1}\times(0,\infty), \medskip\\
u'(x, 0) = u_1(x), u(x, 0) = u_0(x), \ & {\rm in} \ \ \Omega.
\end{array} }} \right.\nonumber
\end{equation}

 When $M(t)=a+b\|\nabla u(t)\|_2^2+\displaystyle\sigma\int_\Omega\nabla u(t)\nabla u_t(t){\rm d}x$, with Balakrishnan-Taylor damping $(\sigma = 0)$, equation \eqref{1.1} describes the motion of a deformable solid with a hereditary effect. The more details about the stability and the well-posedness of the system, we refer the readers to (\cite{ALV2009}, \cite{CAVA2001}, \cite{CAV2008}, \cite{MUS2012}). On the contrary, with Balakrishnan-Taylor damping $(\sigma \neq 0)$ and $g=0$, was initially proposed by Balakrishnan and Taylor in 1989 \cite{BAL1989} and Bass and Zes \cite{BAS1991}. It is used to solve the spillover problem. The related problems also concerned by Tatar and Zarai in (\cite{TA2011}, \cite{ZA2010}, \cite{ZA2013}), Mu in \cite{MU2014} and  Wu in \cite{WU2015}.

Motivated by the above works, we intend to study the local existence, global existence, general decay and blow-up of solutions to problem \eqref{1.1}.
The main difficulties we encounter here arise from the simultaneous appearance of the viscoelastic term, the strong damping term, the Balakrishnan-Taylor damping, the nonlinear source term and the nonlinear boundary damping term. We will show that if the initial data is in the stable set, the solution is global. Otherwise, if the initial data is in the unstable set, by the concavity technique, we get the solution will blow up in a finite time (we can see \cite{LIU2010} for more details). Then, for the relaxation function $g$ satisfies $(G2)$, by introducing suitable energy and Lyapunov functionals, we find that the energy decays at the same rate of $g$, which are not necessarily decaying like polynomial or exponential functions.

The paper is organized as follows. In Section 2, we present some notations and material needed for our work. In Section 3, we establish the local existence. In section 4, the global existence for solutions is given. In Section 5, we prove the general decay result. A finite time blows-up result for initial data in the unstable set is given in Section 6.

\section{Preliminaries }\label{2}
\setcounter{equation}{0}
In this section we present some material needed in the proof of our result. For the relaxation function $g$, we assume

$(G1)$ $g$: $\mathbb{R}_{+}\longrightarrow \mathbb{R}_{+}$ is a nonincreasing differentiable function satisfying
$$g(0)>0,\quad a-\int_{0}^{\infty}g(s){\rm d}s=l>0.$$

$(G2)$ There exists a nonincreasing differentiable function $\xi: \mathbb{R}_{+}\longrightarrow \mathbb{R}_{+}$ such that
$$g'(s)\leq-\xi(s) g(s),\quad \forall s \in\mathbb{R}_{+}.$$

We introduce the set $H_{\Gamma_0}^{1}=\{u\in H^{1}(\Omega)|u_{|_{\Gamma_0}}=0\}$, $\mathcal{V}=H_{\Gamma_{0}}^{1}(\Omega)\cap L^{m}(\Gamma_{1})$. By $(\cdot,\cdot)$ we denote the scalar product in $L^{2}(\Omega)$; i.e.,
$$(u,v)(t)=\int_{\Omega}u(x,t)v(x,t){\rm d}x.$$
Also, by $\|\cdot\|_{q}$ we mean the $L^{q}(\Omega)$ norm for $1\leq q\leq\infty$, by $\|\cdot\|_{q,\Gamma_{1}}$ the $L^{q}(\Gamma_{1})$ norm.
We will also use the embedding (see \cite{ADA1975}):
$H_{\Gamma_0}^{1}\hookrightarrow L^{q}(\Gamma_1)$,\quad $2\leq q\leq \overline q$, where $\overline q$ satifies
\begin{align}
\overline q =\left\{ {{\begin{array}{*{20}l} \frac {2(N-1)}{N-2} > 1, \ & if \ \ N \geq3, \medskip\\ +\infty, \ & if \ \ N=1,2. \end{array} }} \right.\nonumber
\end{align}
We introduce the following functionals as in (\cite{BM2006}, \cite{MT2003}, \cite{MT2007})
\begin{align}\label{2.1}
I(t):=I(u(t))=\left(a-\int_{0}^{t}g(s){\rm d}s\right)\|\nabla u(t)\|_{2}^{2}+(g\circ\nabla u)(t)-\| u(t)\|_{p}^{p}+\frac{b}{2}\|\nabla u(t)\|_2^4,
\end{align}
\begin{align}\label{2.2}
J(t):=J(u(t))=\frac{1}{2}\left[\left(a-\int_{0}^{t}g(s){\rm d}s\right)\|\nabla u(t)\|_{2}^{2}+\frac{b}{2}\|\nabla u(t)\|_2^4+(g\circ\nabla u)(t)\right]-\frac{1}{p}\| u(t)\|_{p}^{p},
\end{align}
\begin{align}\label{2.3}
E(t):=E(u(t),u_{t}(t))=J(t)+\frac{1}{2}\|u_{t}(t)\|_{2}^{2}+\frac{1}{2}\|u_{t}(t)\|_{2,\Gamma_{1}}^{2},
\end{align}
where
$$(g\circ\nabla u)(t)=\int_{0}^{t}g(t-s)\|\nabla u(t)-\nabla u(s)\|_{2}^{2}ds\geq0.$$
We define the potential well depth as
\begin{align}\label{2.4}
d(t)=\inf\limits_{u\in H_{\Gamma_0}^{1}(\Omega)\backslash\{0\}}\sup\limits_{\lambda\geq 0}J(\lambda u).
\end{align}

\section{Local existence result}\label{3}
In this section we will prove the local existence and the uniqueness of the solution of  problem  \eqref{1.1}. The process is closely related to $\left(\cite{GH2008}, {\rm Theorem} 2.1\right)$ that Gerbi and Said-Houair have proved, where no memory term and Balakrishnan-Taylor damping were present.
\begin{theo}\label{th3.1}
 Assume $(G1)$ and $(G2)$ hold, $2\leq p\leq \overline q$ and ${\rm max}\left\{2,\frac{\overline q}{\overline q+1-p}\right\}\leq m\leq \overline q$. Then given $u_0\in H_{\Gamma_0}^{1}(\Omega)$ and $u_{1}\in L^{2}(\Omega)$, there exist $T>0$ and a unique solution $u$ of problem \eqref{1.1} on $(0,T)$ such that
$$u\in C\left([0,T],H_{\Gamma_0}^{1}(\Omega)\right)\cap C^{1}\left([0,T],L^{2}(\Omega)\right),$$
$$u_{t}\in L^2\left(0,T;H_{\Gamma_0}^{1}(\Omega)\right)\cap L^{m}\left((0,T)\times \Gamma_1\right).$$
\end{theo}

By using the Fadeo-Galerkin approxiations and  the well-known contraction mapping theorem we can prove this theorem. In order to define the function for which a fixed point exists, for $u\in C\left([0,T],H_{\Gamma_0}^{1}(\Omega)\right)\cap C^{1}\left([0,T],L^{2}(\Omega)\right)$, we consider the following problem
\begin{equation}\label{3.1}
\left\{ {{\begin{array}{*{20}l} \displaystyle v_{tt}-N(t)\Delta v +\int_{0}^{t}g(t-s)\Delta v(x,s){\rm d}s + \alpha v_{t}-\beta\Delta v_{t} =|u|^{p-2}u , \ & x\in \Omega,  t > 0, \medskip\\
\displaystyle v(x,t) =0, \ & x\in \Gamma_{0},  t > 0, \medskip\\
\displaystyle v_{tt}(x, t) = -N(t)\frac {\partial v}{\partial \nu}(x,t)+\int_{0}^{t}g(t-s)\frac{\partial v}{\partial \nu}(x,s){\rm d}s-\beta\frac {\partial v_{t}}{\partial \nu}(x,t)\medskip\\\quad\quad\quad\quad\ \ -\gamma |v_t|^{m-2}v_{t}(x,t), \ & x \in \Gamma_{1}, t > 0, \medskip\\
v(x, 0) = u_0(x), v_t(x, 0) = u_1(x), \ & x \in  \Omega,
\end{array} }} \right.
\end{equation}
where $N(t)=a+b\|\nabla v(t)\|_2^2+\displaystyle\sigma\int_\Omega\nabla v(t)\nabla v_t(t){\rm d}x$.
\begin{definition}\label{Definition 3.3}
A function $v(x,t)$ such that
$$v\in L^{\infty}\left(0, T; H_{\Gamma_0}^{1}(\Omega)\right),$$
$$v_t\in L^{2}\left(0, T; H_{\Gamma_0}^{1}(\Omega)\right)\cap L^{m}\left((0,T)\times \Gamma_1\right),$$
$$v_t\in L^{\infty}\left(0, T; H_{\Gamma_0}^{1}(\Omega)\right)\cap L^{\infty}\left((0,T)\times L^2{(\Gamma_1)}\right),$$
$$v_{tt}\in L^{\infty}\left(0, T; L^2(\Omega)\right)\cap L^{\infty}\left((0,T)\times L^2{(\Gamma_1)}\right),$$
$$v(x,0)=u_{0}(x),\quad v_{t}(x,0)=u_{1}(x),$$
is a generalized solution to problem \eqref{3.1} if for any function $\omega\in H_{\Gamma_0}^{1}(\Omega)\cap L^{m}(\Gamma_1)$ and $\varphi\in C^{1}(0,T)$ with $\varphi(T)=0$, we get
\begin{align}
\int_{0}^{T}(|u|^{p-2}u,w)(t)\varphi(t){\rm d}t=&\int_{0}^{T}\bigg[(v_{tt},\omega)(t)+a(\nabla v,\nabla \omega)(t)+b\int_\Omega\|\nabla v(t)\|_2^2\Delta v(t) \omega(t){\rm d}x\bigg.\nonumber\\
&\left.+\frac{\sigma}{2}\int_\Omega\frac{{\rm d}}{{\rm d}t}\|\nabla v(t)\|_2^2\Delta v(t) \omega(t){\rm d}x +\alpha(v_t,\omega)(t)+\beta(\nabla v_t,\nabla \omega)(t)\right.\nonumber\\
&\left.-\int_{0}^{t}g(t-s)\nabla v(s)\cdot\nabla \omega(t){\rm d}s\right]\varphi(t){\rm d}t\nonumber\\
&+\int_{0}^{T}\varphi(t)\int_{\Gamma_1}\left[v_{tt}(t)+\gamma|v_t(t)|^{m-2}v_{t}(t)\right]\omega {\rm d}\sigma {\rm d}t.\nonumber
\end{align}
\end{definition}
\begin{lemma}\label{Lemma 3.4}
 Assume $(G1)$, $(G2)$, $2\leq p\leq\overline q$ and $2\leq m\leq\overline q$ hold. Let $u_0\in H^{2}(\Omega)\cap \mathcal{V},$ $u_1\in H^{2}(\Omega)$ and $f\in H^{1}(0,T; L^{2}(\Omega))$, then for any $T>0$, there exists a unique generalized solution $v(t,x)$ of problem \eqref{3.1}.
\end{lemma}
{\bf Proof.} We use Faedo-Galerkin method, some difficulties appear deriving a second-order estimate of $v_{tt}(0)$, inspired by the ideas of Doronin and Larkin in \cite{DL2002} and Cavalcanti et al \cite{CCSM2000}, we set $\widetilde v= v(t,x)-\psi(t,x)$ with $\psi(t,x)=u_0(x)+t u_{1}(x)$. Therefore, we have
\begin{equation}\label{3.3}
\left\{ {{\begin{array}{*{20}l} \displaystyle \widetilde v_{tt}-\widetilde N(t)\Delta \widetilde v +\displaystyle\int_{0}^{t}g(t-s)\Delta \widetilde v(x,s){\rm d}s + \alpha \widetilde v_{t}-\beta\Delta \widetilde v_{t} =f(x,t) \medskip\\ \quad\quad\quad\quad\quad +\widetilde N(t)\Delta\psi-\displaystyle\int_{0}^{t}g(t-s)\Delta\psi(s){\rm d}s-\alpha\psi_t+\beta\Delta\psi_t, \ & x\in \Omega,  t > 0, \medskip\\
\displaystyle \widetilde v(x,t) =0, \ & x\in \Gamma_{0},  t > 0, \medskip\\
\displaystyle \widetilde v_{tt}(x, t) = -\widetilde N(t)\frac {\partial(\widetilde v+\psi)}{\partial \nu}(x,t)+\int_{0}^{t}g(t-s)\frac{\partial (\widetilde v+\psi)}{\partial \nu}(x,s){\rm d}s\medskip\\ \quad\quad\quad\quad\ \ \displaystyle -\beta\frac {\partial (\widetilde v_t+\psi_t)}{\partial \nu}(x,t) -\gamma |\widetilde v_t+\psi_t|^{m-2}(\widetilde v_t+\psi_t)(x,t), \ & x \in \Gamma_{1}, t > 0, \medskip\\
\widetilde v(x, 0) = 0, \widetilde v_t(x, 0) = 0, \ & x \in  \Omega,
\end{array} }} \right.
\end{equation}
where $\widetilde N(t)=a+b\|\nabla \left(\widetilde v+\psi\right)\|_2^2+\displaystyle\frac{\sigma}{2}\frac{{\rm d}}{{\rm d}t}\|\nabla \left(\widetilde v+\psi\right)\|_2^2$. We obtain that if $\widetilde v$ is a solution of problem \eqref{3.3} on $[0,T]$, then $v$ is a solution of problem \eqref{3.1} on $[0,T]$. By using the Faedo-Galerkin method, for every $h\geq1$, let $W_h = {\rm Span}\{\omega_1, \cdot\cdot\cdot ,\omega_h\}$, where $\{\omega_j\}$ is the orthogonal complete system of eigenfunctions of $-\Delta$ in the space $\mathcal{V}$ such that $\|w_j\|_2=1$ for all $j$. Then $\{\omega_j\}$ is orthogonal and complete in $L^{2}(\Omega)\cap L^{2}(\Gamma_1)$. Let
\begin{align}\label{3.4}
\widetilde v_h(t)=\sum\limits_{j=1}^{n}\vartheta_{jh}(t)\omega_j,
\end{align}
where $\widetilde v_h(t)$ are solutions to the finite-dimensional Cauchy problem, we have
\begin{align}\label{3.5}
 &\int_\Omega\widetilde v_{tth}(t)\omega_{j}{\rm d}x+\widetilde N(t)\int_\Omega\nabla(\widetilde v_h+\psi)\cdot\nabla\omega_{j}{\rm d}x
+\alpha\int_\Omega(\widetilde v_h+\psi)_t\omega_{j}{\rm d}x\nonumber \\
 &+\beta\int_\Omega\nabla(\widetilde v_h+\psi)_t\cdot\nabla\omega_{j}{\rm d}x-\int_{0}^{t}g(t-s)\int_\Omega\nabla(\widetilde v_h+\psi)\cdot\nabla\omega_{j}{\rm d}x\nonumber\\
 &+\int_{\Gamma_1}(\widetilde v_{tth}+\gamma|(\widetilde v_h+\psi)_t|^{m-2}(\widetilde v_h+\psi)_t)\omega_j{\rm d}\sigma\nonumber\\
 =&\int_{\Omega}f\omega_j{\rm d}x,
\end{align}
$$\vartheta_{jh}(0)=\vartheta'_{jh}(0)=0,\ \ j=1, \cdot\cdot\cdot ,h.$$
According to \cite{CL1955}, by the Caratheodory theorem, problem \eqref{3.5} has a unique solution $(\vartheta_{jh}(t))\in H^{3}(0,t_h)$ defined on $[0, t_h)$.
To do this, we need two priori estimates to show that for all $h\in N$, $t_h=T$ and those approximations converge to a solution of problem \eqref{3.3}, for the proof of the estimates can be done as in $\left(\cite{GH2008}, {\rm Lemma} 2.2\right)$, we omit the details. Then,
by using Aubin-Lions compactness lemma, the proof can be completed arguing as in \cite{LIO1969}. The uniqueness of this solution follows from the energy inequality, we omit here. This finishes the proof of Lemma \ref{Lemma 3.4}.
Next, we show the following existence result.
\begin{lemma}\label{Lemma 3.2}
Assume $(G1)$, $(G2)$, $2\leq p\leq \overline q$ and ${\rm max}\left\{2,\frac{\overline q}{\overline q+1-p}\right\}\leq m\leq \overline q$ hold. Then given $u_0\in H_{\Gamma_0}^{1}(\Omega)$ and $u_{1}\in L^{2}(\Omega)$, there exist $T>0$ and a unique solution $v$ of problem \eqref{3.1} on $(0,T)$ such that
$$v\in C\left([0,T],H_{\Gamma_0}^{1}(\Omega)\right)\cap C^{1}\left([0,T],L^{2}(\Omega)\right),$$
$$v_{t}\in L^{2}\left(0,T;H_{\Gamma_0}^{1}(\Omega)\right)\cap L^{m}\left((0,T)\times \Gamma_1\right)$$
and satisfies the energy identity
\begin{align}\label{3.50}
&\frac{1}{2}\left[\left(a-\int_{0}^{\tau}g(z){\rm d}z\right)\|\nabla v\|_{2}^{2}+\frac{b}{2}\|\nabla v\|_{2}^{4}+\|v_t\|_{2}^{2}+\|v_t\|_{2,\Gamma_1}^{2}+(g\circ\nabla v)(\tau)\right]_{s}^{t}+\alpha\int_{s}^{t}\|v_t(\tau)\|_{2}^{2}{\rm d}\tau\nonumber\\
&+\frac{\sigma}{4}\int_{s}^{t}\left(\frac{{\rm d}\|\nabla v(z)\|}{{\rm d}z}\right)^{2}{\rm d}z+\beta\int_{s}^{t}\|\nabla v_t(\tau)\|_{2}^{2}{\rm d}s +\gamma\int_{s}^{t}\|v_t(\tau)\|_{m,\Gamma_1}^{m}{\rm d}\tau-\frac{1}{2}\int_{s}^{t}(g'\circ\nabla v)(\tau){\rm d} \tau\nonumber\\
&+\frac{1}{2}\int_{s}^{t}g(\tau)\|\nabla v(\tau)\|_{2}^{2}{\rm d}\tau\nonumber\\
=&\int_{s}^{t}\int_{\Omega}|u(\tau)|^{p-2}u(\tau)v_{t}(\tau){\rm d}\tau {\rm d}x
\end{align}
for $0\leq s\leq t\leq T$.
\end{lemma}
%
%
{\bf Proof.} First, by a sequence $(u^{\iota})_{\iota\in N}\subset C^\infty([0,T]\times\overline\Omega)$ by standard convolution arguments (see \cite{{BRE1983}}), we approximate $u\in C\left([0,T],H_{\Gamma_0}^{1}(\Omega)\right)\cap C^{1}\left([0,T],L^{2}(\Omega)\right)$ endowed with the standard norm
\begin{align}
\|u\|= \max\limits_{t\in [0,T]}\left(\|u_t(t)\|_2+\|u(t)\|_{H^{1}(\Omega)}\right).\nonumber
\end{align}
We obtain $f(u^\iota)=|u^\iota|^{p-2}u^\iota\in H^1(0,T;L^2(\Omega)).$ Next, we approximate the initial data $u_1\in L^2(\Omega)$ by a sequence $\{u_1^\iota\}\in C_0^\infty(\Omega)$. At last, because the space $H^2(\Omega)\cap V \cap H_{\Gamma_0}^1(\Omega)$ is dense in $ H_{\Gamma_0}^1(\Omega)$ for the $H^1$ norm, we approximate $u_0\in H_{\Gamma_0}^1(\Omega)$ by a sequence $\{u_0^\iota\}\subset H^2(\Omega)\cap V \cap H_{\Gamma_0}^1(\Omega)$.

We now consider the following problem
\begin{equation}\label{3.27}
\left\{ {{\begin{array}{*{20}l} \displaystyle v_{tt}^\iota-N^\iota(t)\Delta v^\iota +\int_{0}^{t}g(t-s)\Delta v^\iota(x,s){\rm d}s + \alpha v_{t}^\iota-\beta\Delta v_{t}^\iota =|u^\iota|^{p-2}u^\iota , \ & x\in \Omega,  t > 0, \medskip\\
\displaystyle v^\iota(x,t) =0, \ & x\in \Gamma_{0},  t > 0, \medskip\\
v_{tt}^\iota(x, t) = -\displaystyle N^\iota(t)\frac {\partial v^\iota}{\partial \nu}(x,t)+\displaystyle\int_{0}^{t}g(t-s)\frac{\partial v^\iota}{\partial \nu}(x,s){\rm d}s-\beta\frac {\partial v_{t}^\iota}{\partial \nu}(x,t)\medskip\\\quad\quad\quad\quad\ \ -\gamma |v_t^\iota|^{m-2}v_{t}^\iota(x,t), \ & x \in \Gamma_{1}, t > 0, \medskip\\
v^\iota(x, 0) = u_0^\iota(x), v_t^\iota(x, 0) = u_1^\iota(x), \ & x \in  \Omega,
\end{array} }} \right.
\end{equation}
where $N^\iota(t)=a+b\|\nabla v^\iota(t)\|_2^2+\displaystyle\sigma\int_\Omega\nabla v^\iota(t)\nabla v_t^\iota(t){\rm d}x$. As every hypothesis of Lemma \ref{Lemma 3.4} is satisfied, we can find a sequence of unique solutions $\{v_k\} $ of problem \eqref{3.27}. Next, we shall prove that $\{v^\iota,v_t^\iota \}$ is a Cauchy sequence in the space
\begin{align}
Y_T =\left\{(v,v_t): v\in C\left([0,T], H_{\Gamma_0}^1(\Omega)\right)\cap C^1\left([0,T], L^2(\Omega)\right)\right.,\nonumber\\
\left. v_t\in L^2\left(0,T; H_{\Gamma_0}^1(\Omega)\right)\cap L^m ((0,T)\times\Gamma_1)\right\}\nonumber
\end{align}
endowed with the norm
\begin{align}
\|(v,v_t)\|_{Y_T}^2=\max\limits_{t\in [0,T]}\left[\|v_t\|_2^2+\|\nabla v\|_2^2\right]+\|v_t\|_{L^m\left((0,T)\times\Gamma_1\right)}^2+\int_0^t\|\nabla v_t(s)\|_2^2{\rm d}s. \nonumber
\end{align}
To achieve this, we set $U=u^\iota-u^{\iota'}$, $V=v^\iota-v^{\iota'}$. It's easy to see that $V$ satisfies
\begin{equation}
\left\{ {{\begin{array}{*{20}l} \displaystyle V_{tt}-N_1(t)\Delta V +\int_{0}^{t}g(t-s)\Delta V(x,s){\rm d}s \\\quad\quad\quad\ \ + \alpha V_{t}-\beta\Delta V_{t} =|u^\iota|^{p-2}u^\iota\medskip-|u^{\iota'}|^{p-2}u^{\iota'} , \ & x\in \Omega,  t > 0, \nonumber \medskip\\
\displaystyle V(x,t) =0, \ & x\in \Gamma_{0},  t > 0, \medskip\\
V_{tt}(x, t) = -\displaystyle N_1(t)\frac {\partial V}{\partial \nu}(x,t)+\displaystyle\int_{0}^{t}g(t-s)\frac{\partial V}{\partial \nu}(x,s){\rm d}s-\beta\frac {\partial V_{t}}{\partial \nu}(x,t)\medskip\\\quad\quad\quad\quad\ \ -\gamma \left(|v_t^\iota|^{m-2}v_{t}^\iota(x,t)-|v_t^{\iota'}|^{m-2}v_{t}^{\iota'}(x,t)\right), \ & x \in \Gamma_{1}, t > 0, \medskip\\
V(x, 0) = u_0^\iota-u_0^{\iota'}, V_t(x, 0) = u_1^\iota(x)-u_1^{\iota'}, \ & x \in  \Omega,
\end{array} }} \right.
\end{equation}
where $N_1(t)=a+b\|\nabla V(t)\|_2^2+\displaystyle\sigma\int_\Omega\nabla V(t)\nabla V_t(t){\rm d}x$. By multiplying the above differential equation by $V_t$ and integrating over $\Omega\times(0,t)$, we arrive at
\begin{align}
&\frac{1}{2}\left[\left(a-\int_0^tg(s){\rm d}s\right)\|\nabla V(t)\|_{2}^{2}+\frac{b}{2}\|\nabla V(t)\|_2^4+\|V_{t}(t)\|_{2}^{2}+\|V_{t}(t)\|_{2,\Gamma_1}^{2}\right]
+\alpha\int_{0}^{t}\|V_t(s)\|_2^2{\rm d}s\nonumber\\
&+\beta\int_{0}^{t}\|\nabla V_t(s)\|_2^2{\rm d}s-\int_{0}^{t}(g'\circ\nabla V)(s){\rm d}s+(g\circ\nabla V)(t)+\int_{0}^{t}\int_\Omega g(s)|\nabla V(s)|^2{\rm d}x{\rm d}s\nonumber\\
&+\frac{\sigma}{4}\|\nabla V(t)\|_2^2+\gamma\int_{0}^{t}\int_{\Gamma_1}\left(|v_t^\iota|^{m-2}v_{t}^\iota-|v_t^{\iota'}|^{m-2}v_{t}^{\iota'}\right)\left(v_{t}^\iota-v_{t}^{\iota'}\right){\rm d}\sigma{\rm d}s\nonumber\\
=&\frac{1}{2}\left[\left(a-\int_0^tg(s){\rm d}s\right)\|\nabla V(0)\|_{2}^{2}+\frac{b}{2}\|\nabla V(0)\|_2^4+\|V_{t}(0)\|_{2}^{2}+\|V_{t}(0)\|_{2,\Gamma_1}^{2}\right]\nonumber\\
&+\int_{0}^{t}\int_{\Omega}\left(|u^\iota|^{p-2}u^\iota-|u^{\iota'}|^{p-2}u^{\iota'}\right)\left(v_t^\iota-v_t^{\iota'}\right){\rm d}x{\rm d}s, \ \ \forall t\in (0,T).\nonumber
\end{align}
By the algebraic inequality
\begin{align}\label{3.271}
\forall m\geq2, \ \ \exists c>0, \ \ \forall a,b\in \mathbb{R}, \ \ \left(|a|^{m-2}a-|b|^{m-2}b\right)(a-b)\geq c|a-b|^m,
\end{align}
we get
\begin{align}
&\frac{1}{2}\left[\left(a-\int_0^tg(s){\rm d}s\right)\|\nabla V(t)\|_{2}^{2}+\frac{b}{2}\|\nabla V(t)\|_2^4+\|V_{t}(t)\|_{2}^{2}+\|V_{t}(t)\|_{2,\Gamma_1}^{2}\right]
+\alpha\int_{0}^{t}\|V_t(s)\|_2^2{\rm d}s\nonumber\\
&+c_1\|V_t\|_{m,\Gamma_1}^m+\displaystyle \frac{\sigma}{4}\|\nabla V(t)\|_2^2+\beta\int_{0}^{t}\|\nabla V_t(s)\|_2^2{\rm d}s-\int_{0}^{t}(g'\circ\nabla V)(s){\rm d}s+(g\circ\nabla V)(t)\nonumber\\
&+\int_{0}^{t}\int_\Omega g(s)|\nabla V(s)|^2{\rm d}x{\rm d}s\nonumber\\
\leq&\frac{1}{2}\left[\left(a-\int_0^tg(s){\rm d}s\right)\|\nabla V(0)\|_{2}^{2}+\frac{b}{2}\|\nabla V(0)\|_2^4+\|V_{t}(0)\|_{2}^{2}+\|V_{t}(0)\|_{2,\Gamma_1}^{2}\right]\nonumber\\
&+\int_{0}^{t}\int_{\Omega}\left(|u^\iota|^{p-2}u^\iota-|u^{\iota'}|^{p-2}u^{\iota'}\right)\left(v_t^\iota-v_t^{\iota'}\right){\rm d}x{\rm d}s, \ \ \forall t\in (0,T).\nonumber
\end{align}
Then, by using Young's inequality, Gronwall inequality and the result of Georgiev and Todorova \cite{GT1994}, there exists a constant $C$ depending only on $\Omega$ and $p$ such that
\begin{align}
\|V\|_{Y_T}\leq C\left(\|\nabla V(0)\|_{2}^{2}+\|V_{t}(0)\|_{2}^{2}+\|V_{t}(0)\|_{2,\Gamma_1}^{2}\right)+CT\|U\|_{Y_T}.\nonumber
\end{align}
By the above notation, we can obtain $V(0)=u_0^\iota-u_0^{\iota'}$, $V_t(0)=u_1^\iota-u^{\iota'}$. Thus, we conclude that $\{\left(v^\iota,v_t^\iota\right)\}$ is a Cauchy sequence in $Y_T$, because $\{u_0^\iota\}$ is a converging sequence in $H_{\Gamma_0}^1(\Omega)$,  $\{u_1^\iota\}$ is a converging sequence in $L^2(\Omega)$ and  $\{u^\iota\}$ is a converging sequence in $C\left([0,T], H_{\Gamma_0}^1(\Omega)\right)\cap C^1\left([0,T], L^2(\Omega)\right)$.
So, we conclude that $\{\left(v^\iota,v_t^\iota\right)\}$ converges to a limit $(v,v_t)\in Y_T$. As done by \cite{GT1994}, we can prove that this limit is a weak solution of problem \eqref{3.1}. This completes the proof of Lemma \ref{Lemma 3.2}.

Now we are ready to prove the local existence result.

{\bf Proof of Theorem \ref{th3.1}.}  For $R>0$ large and $T>0$, we define the convex closed subset of $Y_T$, we define a class of functions as
\begin{align}
X_T=\{(v,v_t)\in Y_T : v(0)=u_0, v_t(0)=u_1 \ \ {\rm and}\ \  \|v\|_{Y_T}\leq R\}.\nonumber
\end{align}
Then, Lemma \ref{Lemma 3.2} implies that, for any $u\in X_T$, we may define $v=\Phi (u)$ to be a unique solution of \eqref{3.1} corresponding to $u$. We would like to show that $\Phi$ is a contraction map satisfying $\Phi(X_T)\subset X_T$ for a suitable $T>0$.
Firstly, we show $\Phi(X_T)\subset X_T$. For this, by using the energy identity, we have
\begin{align}\label{3.28}
&l\|\nabla v(t)\|_2^2+\|v_t(t)\|_2^2+\|v_t(t)\|_{2,\Gamma_1}^2+2\alpha\int_0^t\|v_t(s)\|_2^2{\rm d}s+2\beta\int_0^t\|\nabla v_t(s)\|_2^2{\rm d}s\nonumber\\
&+2\gamma\int_0^t\|v_t(s)\|_{m,\Gamma_1}^m{\rm d}s\nonumber\\
\leq& \|\nabla u_0\|_2^2+\|u_1\|_2^2+\|u_1\|_{2,\Gamma_1}^2+2\int_0^t\int_\Omega|u(s)|^{p-2}u(s)v_t(s){\rm d}x{\rm d}s.
\end{align}
By the Holder's inequality, we estimate the last term in the right-hand side of the inequality \eqref{3.28} as follows
\begin{align}
&\int_0^t\int_\Omega|u(s)|^{p-2}u(s)v_t(s){\rm d}x{\rm d}s\nonumber\\
\leq& \int_0^t\|u(s)\|_{2N/(N-2)}^{p-1}\|v_t(s)\|_{2N/\left(3N-Np+2(p-1)\right)} {\rm d}s.\nonumber
\end{align}
As $p\leq\frac{2N}{N-2}$, we get
\begin{align}
\frac{2N}{\left(3N-Np+2(p-1)\right)}\leq\frac{2N}{N-2}.\nonumber
\end{align}
Thus, by Young's and Sobolev's inequalities, for all $\delta>0, t\in(0,T)$, there exists a positive constant $C(\delta)$ such that
\begin{align}
\int_0^t\int_\Omega|u(s)|^{p-2}u(s)v_t(s){\rm d}x{\rm d}s
\leq C(\delta)tR^{2(p-1)}+\delta\int_0^t\|\nabla v_t(s)\|_2^2{\rm d}s.\nonumber
\end{align}
Inserting the last estimate in the inequality \eqref{3.28} and choosing $\delta$ small enough, we have
\begin{align}
\|v\|_{Y_T}^2\leq\frac{1}{2}R^2+CTR^{2(p-1)}.\nonumber
\end{align}
Then, choose $T$ small enough, we get $\|v\|_{Y_T}\leq R$, which shows that $\Phi$ maps $X_T$ into itself.

Next, we verify that $\Phi$ is a contraction. To this end, we set $U=u-\overline u$ and $V=v-\overline v$, where $v=\Phi(u)$ and $\overline v=\Phi(\overline u)$.
It's straightforward to verify that $V$ satisfies
\begin{equation}\label{3.29}
\left\{ {{\begin{array}{*{20}l} \displaystyle V_{tt}-N_2(t)\Delta V +\int_{0}^{t}g(t-s)\Delta V(x,s){\rm d}s \\\quad\quad\ \ \ \ + \alpha V_{t}-\beta\Delta V_{t} =|u|^{p-2}u\medskip-|\overline u|^{p-2}\overline u , \ & x\in \Omega,  t > 0,  \medskip\\
\displaystyle V(x,t) =0, \ & x\in \Gamma_{0},  t > 0, \medskip\\
V_{tt}(x, t) = -\displaystyle N_2(t)\frac {\partial V}{\partial \nu}(x,t)+\displaystyle\int_{0}^{t}g(t-s)\frac{\partial V}{\partial \nu}(x,s){\rm d}s-\beta\frac {\partial V_{t}}{\partial \nu}(x,t)\medskip\\\quad\quad\quad\quad\ \ -\gamma \left(|v_t|^{m-2}v_{t}(x,t)-|\overline {v_t}|^{m-2}\overline {v_{t}}(x,t)\right), \ & x \in \Gamma_{1}, t > 0, \medskip\\
V(x, 0) = 0, V_t(x, 0) = 0, \ & x \in  \Omega,
\end{array} }} \right.
\end{equation}
where $N_2(t)=a+b\|\nabla V(t)\|_2^2+\displaystyle\sigma\int_\Omega\nabla V(t)\nabla V_t(t){\rm d}x$. By multiplying the differential equation in \eqref{3.29} by $V_t$ and integrating over $\Omega\times (0,t)$, we arrive at
\begin{align}
&\frac{1}{2}\left[\left(a-\int_0^tg(s){\rm d}s\right)\|\nabla V(t)\|_{2}^{2}+\frac{b}{2}\|\nabla V(t)\|_2^4+\|V_{t}(t)\|_{2}^{2}+\|V_{t}(t)\|_{2,\Gamma_1}^{2}\right]
+\alpha\int_{0}^{t}\|V_t(s)\|_2^2{\rm d}s\nonumber\\
&+\beta\int_{0}^{t}\|\nabla V_t(s)\|_2^2{\rm d}s-\frac{1}{2}\int_{0}^{t}(g'\circ\nabla V)(s){\rm d}s+\frac{1}{2}(g\circ\nabla V)(t)+\frac{1}{2}\int_{0}^{t}\int_\Omega g(s)|\nabla V(s)|^2{\rm d}x{\rm d}s\nonumber\\
&+\displaystyle \frac{\sigma}{4}\|\nabla V(t)\|_2^2+\gamma\int_{0}^{t}\int_{\Gamma_1}\left(|v_t|^{m-2}v_{t}-|\overline {v_t}|^{m-2}\overline {v_{t}}\right)\left(v_{t}-\overline {v_{t}}\right){\rm d}\sigma{\rm d}s\nonumber\\
=&\int_{0}^{t}\int_{\Omega}\left(|u|^{p-2}u-|\overline u|^{p-2}\overline u\right)\left(v_t-\overline {v_t}\right){\rm d}x{\rm d}s, \ \ \forall t\in (0,T).\nonumber
\end{align}
Similar to the discussion in \cite{GH2008}, we can get
\begin{align}\label{3.35}
\|V\|_{Y_T}^2\leq CR^{p-2}T^{1/2}\|U\|_{Y_T}^2.
\end{align}
By choosing $T$ so small that $CR^{p-2}T^{1/2}<1$, estimate \eqref{3.35} shows that $\Phi$ is a contraction. The contraction mapping theorem then guarantees the existence of a unique $v$ satisfying $v=\Phi(v)$. The proof of Theorem \ref{th3.1} is now completed.

\section{Global existence} \label{4}
\setcounter{equation}{0}
In this section, for the initial data in the stable set, we show the solution is global. We need the following definition and lemmas.
\begin{definition}\label{Definition 4.1}
Let $2\leq p\leq \overline q$ and ${\rm max}\left\{2,\frac{\overline q}{\overline q+1-p}\right\}\leq m\leq \overline q$. 
We denote $u$ as the solution of  \eqref{1.1}. We define
\begin{align}
T_{{\rm max}}= \sup\{T>0, u=u(t) \ \ {\rm exists}\ \ {\rm on}\ \ [0,T]\}.\nonumber
\end {align}
\end{definition}
\begin{lemma}\label{Lemma 4.1}
Suppose that $(G1)$, $(G2)$, $2\leq p\leq \overline q$ and ${\rm max}\left\{2,\frac{\overline q}{\overline q+1-p}\right\}\leq m\leq \overline q$ hold. If $u$ is the solution of \eqref{1.1}, then
\begin{align}\label{4.1}
E'(t)=&\frac{1}{2}(g'\circ\nabla u)(t)-\frac{1}{2}g(t)\|\nabla u(t)\|_{2}^{2}-\sigma\left(\frac{1}{2}\frac{\rm d}{\rm dt}\|\nabla u(t)\|_{2}^{2}\right)^2-\alpha\|u_{t}(t)\|_{2}^{2}\nonumber\\
-&\beta\|\nabla u_{t}(t)\|_{2}^{2}-\gamma\| u_{t}(t)\|_{m,\Gamma_{1}}^{m}\leq0,
\end {align}
for almost every $t\in [0,T_{{\rm max}}).$
\end{lemma}
{\bf Proof.}  In view of $(G1)$, multiplying the differential equation in \eqref{1.1} by $u_t$ and integrating by parts over $\Omega$, we obtain the result.
\begin{lemma}\label{Lemma 4.2}
Let $(u_0,u_1)\in H_{\Gamma_0}^1\times L^{2}(\Omega)$ be given, suppose that $(G1)$ and $2\leq p\leq \overline q$ hold, such that
\begin{equation}\label{4.2}
\left\{ {{\begin{array}{*{20}l} \displaystyle \frac{C_{\ast}^{p}}{l}\left(\frac{2p}{l(p-2)}E(0)\right)^{\frac{p-2}{2}}<1,\medskip\\
\displaystyle I(u_0)>0,
\end{array} }} \right.
\end{equation}
where $C_\ast$ is the best poincare constant, then $I(u(t))>0$, $\forall t\in [0,T_{{\rm max}})$.
\end{lemma}
{\bf Proof.} Since $I(0)>0$, there exists a $T_{\ast}<T_{{\rm max}}$ such that $I(t)>0$ for all $t\in [0,T_{\ast})$. This implies that
\begin{align}\label{4.3}
J(t)=&\frac{p-2}{2p}\left[\left(a-\int_{0}^{t}g(s){\rm d}s\right)\|\nabla u(t)\|_{2}^{2}+(g\circ\nabla u)(t)+\frac{b}{2}\|\nabla u(t)\|_{2}^{4}\right]+\frac{1}{p}I(t)\nonumber\\
\geq&\frac{p-2}{2p}\left[\left(a-\int_{0}^{t}g(s){\rm d}s\right)\|\nabla u(t)\|_{2}^{2}+(g\circ\nabla u)(t)+\frac{b}{2}\|\nabla u(t)\|_{2}^{4}\right].
\end {align}
Thus, by $(G1)$, \eqref{2.1}, \eqref{4.1}, and \eqref{4.3}, we have
\begin{align}\label{4.4}
l\|\nabla u\|_{2}^{2}\leq\left(a-\int_{0}^{t}g(s){\rm d}s\right)\|\nabla u(t)\|_{2}^{2}\leq\frac{2p}{p-2}J(t)\leq\frac{2p}{p-2}E(t)\leq\frac{2p}{p-2}E(0),
\end {align}
for all $t\in [0,T_\ast)$. This combines with the Soblev imbedding, $(G1)$ and \eqref{4.2}, implies that
\begin{align}
\|u\|_{p}^{p}&\leq C_{\ast}^{p}\|\nabla u\|_{2}^{p}\leq\frac{C_{\ast}^{p}}{l}\|\nabla u\|_{2}^{p-2}l\|\nabla u\|_{2}^{2}\nonumber \\
&\leq\frac{C_{\ast}^{p}}{l}\left(\frac{2p}{l(p-2)}E(0)\right)^{\frac{p-2}{2}}l\|\nabla u\|_{2}^{2}<\left(a-\int_{0}^{t}g(s){\rm d}s\right)\|\nabla u\|_{2}^{2},\nonumber
\end {align}
for all $t\in [0,T_\ast)$. Therefore, $I(t)>0$, for all $t\in [0,T_\ast)$. By repeating this procedure and using the fact that
\begin{align}
\lim\limits_{t\rightarrow T_\ast}\frac{C_{\ast}^{p}}{l}\left(\frac{2p}{l(p-2)}E(0)\right)^{\frac{p-2}{2}}<1,\nonumber
\end {align}
where $T_{\ast}$ is extended to $T_{{\rm max}}$.
\begin{theo}\label{th4.4} Let $(u_0,u_1)\in H_{\Gamma_0}^1\times L^{2}(\Omega)$ be given and satisfying \eqref{4.2}, suppose that $(G1)$, $(G2)$, $2\leq p\leq \overline q$ and ${\rm max}\left\{2,\frac{\overline q}{\overline q+1-p}\right\}\leq m\leq \overline q$ hold. Then the solution is global and bounded.
\end{theo}
{\bf Proof.}  We have just to check that $$\|\nabla u(t)\|_{2}^{2}+\| u_{t}(t)\|_{2}^{2}$$ is bounded and independently of $t$. To achieve this, by \eqref{2.1}, \eqref{4.1} and \eqref{4.3},
we have
\begin{align}\label{4.5}
E(0)&\geq E(t)=J(t)+\frac{1}{2}\|u_{t}(t)\|_{2}^{2}+\frac{1}{2}\|u_{t}(t)\|_{2,\Gamma_1}^{2}\nonumber \\
&\geq \frac{p-2}{2p}\left[l\|\nabla u(x,t)\|_{2}^{2}+(g\circ \nabla u)(x,t)+\frac{b}{2}\|\nabla u(t)\|_{2}^{4}\right]+\frac{1}{2}\|u_{t}(t)\|_{2}^{2}+\frac{1}{2}\|u_{t}(t)\|_{2,{\Gamma_1}}^{2}+\frac{1}{p}I(t)\nonumber \\&
\geq \frac{p-2}{2p}l\|\nabla u(t)\|_{2}^{2}+\frac{1}{2}\|u_{t}(t)\|_{2}^{2}+\frac{1}{2}\|u_{t}(t)\|_{2,\Gamma_1}^{2},
\end {align}
since $I(t)$ and $(g\circ \nabla u)(t)$ are positive. Therefore
$$\|\nabla u(t)\|_{2}^{2}+\| u_{t}(t)\|_{2}^{2}\leq CE(0),\quad \forall t\in [0,T_{{\rm max}}),$$
where $C$ is a positive constant, which depends only on $p$ and $l$. Then by the definition of $T_{{\rm max}}$, the solution is global, that is $T_{{\rm max}}=\infty$.
\section{Decay of solutions} \label{5}

In this section, we state and prove the general decay result. For positive constants $\varepsilon_1$ and $\varepsilon_2$, we use the following modified functional
\begin{align}\label{5.1}
L(t)=E(t)+\varepsilon_1 G(t)+\varepsilon_2 H(t),
\end {align}
where
\begin{align}\label{5.2}
G(t)=\int_\Omega u_tu{\rm d}x+\int_{\Gamma_1} u_tu{\rm d}\sigma+\frac{\alpha}{2}\|u\|_2^2+\frac{\beta}{2}\|\nabla u\|_2^2+\frac{\sigma}{4}\|\nabla u\|_{2}^{4}
\end {align}
and
\begin{align}\label{5.3}
H(t)=-\int_\Omega u_t\int_0^{t}g(t-s)\left(u(t)-u(s)\right){\rm d}s{\rm d}x-\int_{\Gamma_1}u_t\int_0^{t}g(t-s)\left(u(t)-u(s)\right){\rm d}s{\rm d}\sigma.
\end {align}
It is easy to check that, by using Poincare's inequality, trace inequality, \eqref{2.2}, \eqref{2.3} and for $\varepsilon_1$, $\varepsilon_2$ small enough,
\begin{align}\label{5.4}
\alpha_1L(t)\leq E(t)\leq \alpha_2L(t)
\end {align}
holds for two positive constants $\alpha_1$ and $\alpha_2$.
\begin{theo}\label{th5.1}
Let $(u_0,u_1)\in H_{\Gamma_0}^1\times L^{2}(\Omega)$ be given. Assume that $g$ and $\xi$ satisfy $(G1)$ and $(G2)$. Then, for each $t_0>0$, there exist positive constants $K$ and $k$ such that the solution \eqref{1.1} satisfies, for all $t\geq t_0$,
\begin{align}\label{5.5}
E(t)\leq Ke^{-k\int_{t_0}^t\xi(s){\rm d}s}.
\end {align}
\end{theo}
\begin{lemma}\label{Lemma 5.2}
Under the conditions of Theorem \ref{th5.1}, the functional $G(t)$ defined by \eqref{5.2} satisfies
\begin{align}\label{5.6}
G'(t)\leq &\| u_t(t)\|_2^2+\|u_{t}(t)\|_{2,\Gamma_1}^{2}-\left(\frac{l}{2}-\frac{C\gamma\delta^{-m}}{m}\left(\frac{2pE(0)}{l(p-2)}\right)^{\frac{m-2}{2}}\right)\|\nabla u(t)\|_2^2+\|u(t)\|_p^p\nonumber\\
&+\frac{\gamma(m-1)}{m}\delta^{\frac{m}{m-1}}\|u_t(t)\|_{m,\Gamma_1}^{m}+\frac{a-l}{2l}\left(g\circ\nabla u\right)(t),
\end {align}
for some $\delta>0$.
\end{lemma}
{\bf Proof.} By using the differential equation in \eqref{1.1}, we get
\begin{align}\label{5.7}
G'(t)=&\| u_t(t)\|_{2}^{2}+\int_\Omega u_{tt}(t)u(t){\rm d}x+\int_{\Gamma_1} u_{tt}(t)u(t){\rm d}\sigma+\|u_{t}(t)\|_{2,\Gamma_1}^{2}+\alpha\int_\Omega u_t(t)u(t){\rm d}x\nonumber\\
&+\beta\int_\Omega \nabla u_t(t)\cdot\nabla u(t){\rm d}x+\frac{\sigma}{2}\left(\frac{\rm d}{\rm dt}\|\nabla u(t)\|_2^2\right)\|\nabla u(t)\|_2^2\nonumber\\
=&\| u_t(t)\|_{2}^{2}+\|u_{t}(t)\|_{2,\Gamma_1}^{2}+\int_{\Gamma_1}\frac{\partial u(t)}{\partial \nu}u(t){\rm d}\sigma-\left(a+b\|\nabla u(t)\|_2^2\right)\| \nabla u(t)\|_{2}^{2}\nonumber\\
&-\int_{0}^{t}g(t-s)\int_{\Gamma_1}\frac{\partial u(s)}{\partial \nu}u(s){\rm d}\sigma{\rm d}s+\int_{0}^{t}g(t-s)\int_\Omega\nabla u(s)\cdot\nabla u(t){\rm d}x{\rm d}s\nonumber\\
&+\beta\int_{\Gamma_1}\frac{\partial u_t(t)}{\partial \nu}u(t){\rm d}\sigma+\int_{\Gamma_1} u_{tt}(t)u(t){\rm d}\sigma +\|u(t)\|_p^p\nonumber\\
=&\| u_t(t)\|_2^2+\|u_{t}(t)\|_{2,\Gamma_1}^{2}-\left(a+b\|\nabla u(t)\|_2^2\right)\|\nabla u(t)\|_2^2+\int_\Omega\nabla u(t)\cdot\int_0^tg(t-s)\nabla u(s){\rm d}s{\rm d}x\nonumber\\
&-\gamma\int_{\Gamma_1}|u_t(t)|^{m-2}u_t(t)u(t){\rm d}\sigma+\|u(t)\|_p^p.
\end {align}
We now estimate the right hand side of \eqref{5.7}. For two positive constants $\delta$ and $\eta$, we have the estimates as follows
\begin{align}\label{5.8}
\int_\Omega\nabla u(t)\cdot\int_0^tg(t-s)\nabla u(s){\rm d}s{\rm d}x\leq(\eta+a-l)\|\nabla u(t)\|_2^2+ \frac{a-l}{4\eta}(g\circ \nabla u)(t).
\end {align}
By Young's inequality, trace inequality and combining \eqref{4.4}, we have
\begin{align}\label{5.9}
\left|\int_{\Gamma_1}|u_t(t)|^{m-2}u_t(t)u(t){\rm d}\sigma\right|\leq \frac{\delta^{-m}}{m}\|u(t)\|_{m,\Gamma_1}^{m}+\frac{m-1}{m}\delta^{\frac{m}{m-1}}\|u_t(t)\|_{m,\Gamma_1}^{m},
\end {align}
\begin{align}\label{5.10}
\|u(t)\|_{m,\Gamma_1}^{m}\leq C\|\nabla u(t)\|_2^m\leq C\left(\frac{2pE(0)}{l(p-2)}\right)^{\frac{m-2}{2}}\|\nabla u(t)\|_2^2.
\end {align}
Combining \eqref{5.7}-\eqref{5.10}, we obtain
\begin{align}
G'(t)\leq &\| u_t(t)\|_2^2+\|u_{t}(t)\|_{2,\Gamma_1}^{2}-\left(l-\eta-\frac{C\gamma\delta^{-m}}{m}\left(\frac{2pE(0)}{l(p-2)}\right)^{\frac{m-2}{2}}\right)\|\nabla u(t)\|_2^2+\|u(t)\|_p^p\nonumber\\
&+\frac{\gamma(m-1)}{m}\delta^{\frac{m}{m-1}}\|u_t(t)\|_{m,\Gamma_1}^{m}+\frac{a-l}{4\eta}\left(g\circ\nabla u\right)(t)\nonumber,
\end {align}
Letting $\eta=\displaystyle\frac{l}{2}>0$ in above inequality, we arrive at \eqref{5.6}.
\begin{lemma}\label{Lemma 5.3}
Under the conditions of Theorem \ref{th5.1}, the functional $H(t)$ defined by \eqref{5.3} satisfies
\begin{align}\label{5.111}
H'(t)\leq&\left(\alpha+\delta-\int_0^{t}g(s){\rm d}s\right)\|u_t(t)\|_2^2+\left[\delta+2\delta(a-l)^2+\delta C_{\ast}^{2(p-1)}\left(\frac{2pE(0)}{l(p-2)}\right)^{p-2}\right]\|\nabla u(t)\|_2^2\nonumber\\
&+\left[\frac{a-l}{4\delta}\left(a+\frac{2pbE(0)}{(p-2)l }\right)^2+\frac{a-l}{4l}+\frac{a-l}{2\delta\lambda_1}+\frac{\beta(a-l)}{2}+\left(2\delta+\frac{1}{4\delta}\right)(a-l)^2\right.\nonumber\\
&\left.+\frac{\gamma C_{\ast}^{m}c}{4\delta m}\left(\frac{8pE(0)}{l(p-2)}\right)^{\frac{m-2}{2}}\right](g\circ\nabla u)(t)+\frac{\beta}{2}\|\nabla u_t(t)\|_2^2+(\delta+a-l)\|u_t(t)\|_{2,\Gamma_1}^2\nonumber\\
&+\frac{\gamma\delta(m-1)}{m}(a-l)\|u_t(t)\|_{m,\Gamma_1}^{m}+\frac{a-l}{2\delta\lambda_1}(-g'\circ\nabla u)(t)+\frac{\sigma(p-2)}{2p}E(0)E'(t),
\end {align}
for some $\delta>0$.
\end{lemma}
{\bf Proof.} By using the differential equation in \eqref{1.1}, we get
\begin{align}\label{5.11}
H'(t)=&-\int_\Omega u_{tt}(t)\int_0^{t}g(t-s)\left(u(t)-u(s)\right){\rm d}s{\rm d}x-\int_\Omega u_{t}(t)\int_0^{t}g'(t-s)\left(u(t)-u(s)\right){\rm d}s{\rm d}x\nonumber\\
&-\left(\int_0^{t}g(s){\rm d}s\right)\int_\Omega|u_t(t)|^{2}{\rm d}x-\int_{\Gamma_1}u_{tt}(t)\int_0^{t}g(t-s)\left(u(t)-u(s)\right){\rm d}s{\rm d}\sigma\nonumber\\
&-\int_{\Gamma_1} u_{t}(t)\int_0^{t}g'(t-s)\left(u(t)-u(s)\right){\rm d}s{\rm d}\sigma-\left(\int_0^{t}g(s){\rm d}s\right)\int_{\Gamma_1}|u_t(t)|^{2}{\rm d}\sigma
\nonumber\\
=&-\int_\Omega\left(M(t)\Delta u(t)-\int_0^{t}g(t-s)\Delta u(s){\rm d}s-\alpha u_t(t)+\beta\Delta u_t(t)+|u(t)|^{p-2}u(t)\right)\nonumber\\
&\times\int_0^{t}g(t-s)\left(u(t)-u(s)\right){\rm d}s{\rm d}x-\int_\Omega u_{t}(t)\int_0^{t}g'(t-s)\left(u(t)-u(s)\right){\rm d}s{\rm d}x\nonumber\\
&-\left(\int_0^{t}g(s){\rm d}s\right)\int_\Omega|u_t(t)|^{2}{\rm d}x-\int_{\Gamma_1} u_{t}(t)\int_0^{t}g'(t-s)\left(u(t)-u(s)\right){\rm d}s{\rm d}\sigma\nonumber\\
&-\int_{\Gamma_1}\left(-M(t)\frac {\partial u}{\partial \nu}(x,t)+\int_{0}^{t}g(t-s)\frac{\partial u}{\partial \nu}(x,s){\rm d}s-\beta\frac {\partial u_{t}}{\partial \nu}(x,t)-\gamma |u_t|^{m-2}u_{t}(x,t)\right)\nonumber\\
&\times\int_0^{t}g(t-s)\left(u(t)-u(s)\right){\rm d}s{\rm d}\sigma-\left(\int_0^{t}g(s){\rm d}s\right)\int_{\Gamma_1}|u_t(t)|^{2}{\rm d}\sigma\nonumber\\
=&\int_\Omega\left(a+b\|\nabla u(t)\|_2^2\right)\nabla u(t)\cdot\int_0^{t}g(t-s)\left(\nabla u(t)-\nabla u(s)\right){\rm d}s{\rm d}x\nonumber\\
&+\int_\Omega\sigma\left(\nabla u(t)\nabla u_t(t){\rm d} x\right)\nabla u(t)\cdot\int_0^{t}g(t-s)\left(\nabla u(t)-\nabla u(s)\right){\rm d}s{\rm d}x\nonumber\\
&+\beta\int_\Omega\nabla u_t(t)\cdot\int_0^{t}g(t-s)\left(\nabla u(t)-\nabla u(s)\right){\rm d}s{\rm d}x\nonumber\\
&-\int_\Omega\int_0^{t}g(t-s)\nabla u(s){\rm d}s\int_0^{t}g(t-s)\left(\nabla u(t)-\nabla u(s)\right){\rm d}s{\rm d}x\nonumber\\
&+\alpha\int_\Omega u_t(t)\int_0^{t}g(t-s)\left(u(t)-u(s)\right){\rm d}s{\rm d}x\nonumber\\
&-\int_\Omega|u(t)|^{p-2}u(t)\int_0^{t}g(t-s)\left(u(t)-u(s)\right){\rm d}s{\rm d}x\nonumber\\
&-\int_\Omega u_t(t)\int_0^{t}g'(t-s)\left(u(t)-u(s)\right){\rm d}s{\rm d}x-\left(\int_0^{t}g(s){\rm d}s\right)\|u_t(t)\|_2^2\nonumber\\
&+\gamma\int_{\Gamma_1}|u_t(t)|^{m-2}u_t(t)\int_0^{t}g(t-s)\left(u(t)-u(s)\right){\rm d}s{\rm d}\sigma\nonumber\\
&-\int_{\Gamma_1}u_t(t)\int_0^{t}g'(t-s)\left(u(t)-u(s)\right){\rm d}s{\rm d}\sigma\nonumber\\
&-\left(\int_0^{t}g(s){\rm d}s\right)\int_{\Gamma_1}|u_t(t)|^{2}{\rm d}\sigma\nonumber\\
=&-\left(\int_0^{t}g(s){\rm d}s\right)\|u_t(t)\|_2^2+\sum_{i=1}^{10} I_i.
\end {align}
We now estimate $I_i$, $i=1,\cdot\cdot\cdot ,10$ in the right side of \eqref{5.11}. For $\delta>0$, similar as in \cite{BM2006}, exploiting $E'(t)\leq\displaystyle -\sigma\left(\frac{1}{2}\frac{\rm d}{{\rm d}t}\|\nabla u(t)\|_2^2\right)$ by \eqref{4.1}, \eqref{4.4}, Young's inequality, Hoider's inequality and  Cauchy-Schwarzs's inequality, we get
\begin{align}\label{5.12}
I_1=&\int_\Omega\left(a+b\|\nabla u(t)\|_2^2\right)\nabla u(t)\cdot\int_0^{t}g(t-s)\left(\nabla u(t)-\nabla u(s)\right){\rm d}s{\rm d}\sigma\nonumber\\
\leq&\delta\|\nabla u(t)\|_2^2+\frac{a-l}{4\delta}\left(a+\frac{2pbE(0)}{(p-2)l}\right)^2(g\circ\nabla u)(t),
\end {align}
\begin{align}\label{5.133}
I_2=&\int_\Omega\sigma\left(\nabla u(t)\nabla u_t(t){\rm d} x\right)\nabla u(t)\cdot\int_0^{t}g(t-s)\left(\nabla u(t)-\nabla u(s)\right){\rm d}s{\rm d}x\nonumber\\
\leq&\delta^2\left(\nabla u(t)\nabla u_t(t){\rm d} x\right)^2l\|\nabla u(t)\|_2^2+\frac{1}{4l}\int_\Omega\left(\int_0^{t}g(t-s)\left(\nabla u(t)-\nabla u(s)\right){\rm d}s\right)^2{\rm d}x\nonumber\\
\leq&\frac{-\sigma(p-2)}{2p}E(0)E'(t)+\frac{a-l}{4l}(g\circ\nabla u)(t),
\end {align}
\begin{align}\label{5.13}
I_3=&\beta\int_\Omega\nabla u_t(t)\cdot\int_0^{t}g(t-s)\left(\nabla u(t)-\nabla u(s)\right){\rm d}s{\rm d}x\nonumber\\
\leq&\frac{\beta}{2}\|\nabla u_t(t)\|_2^2+\frac{\beta(a-l)}{2}(g\circ\nabla u)(t),
\end {align}
\begin{align}\label{5.14}
I_4=&-\int_\Omega\int_0^{t}g(t-s)\nabla u(s){\rm d}s\int_0^{t}g(t-s)\left(\nabla u(t)-\nabla u(s)\right){\rm d}s{\rm d}x\nonumber\\
\leq&2\delta(a-l)^2\|\nabla u(t)\|_2^2+\left(2\delta+\frac{1}{4\delta}\right)(a-l)^2(g\circ\nabla u)(t),
\end {align}
\begin{align}\label{5.15}
I_5=&\alpha\int_\Omega u_t(t)\int_0^{t}g(t-s)\left(u(t)-u(s)\right){\rm d}s{\rm d}x\nonumber\\
\leq&\alpha\|u_t(t)\|_2^2+\frac{a-l}{4\delta\lambda_1}(g\circ\nabla u)(t),
\end {align}
\begin{align}\label{5.16}
I_6=&-\int_\Omega|u(t)|^{p-2}u(t)\int_0^{t}g(t-s)\left(u(t)-u(s)\right){\rm d}s{\rm d}x\nonumber\\
\leq&\delta C_{\ast}^{2(p-1)}\left(\frac{2pE(0)}{l(p-2)}\right)^{p-2}\|\nabla u(t)\|_2^2+\frac{a-l}{4\delta\lambda_1}(g\circ\nabla u)(t),
\end {align}
\begin{align}\label{5.17}
I_7=&-\int_{\Gamma_1}u_t(t)\int_0^{t}g'(t-s)\left(u(t)-u(s)\right){\rm d}s{\rm d}\sigma\nonumber\\
\leq&\delta\|u_t(t)\|_2^2+\frac{a-l}{4\delta\lambda_1}(-g'\circ\nabla u)(t),
\end {align}
\begin{align}\label{5.18}
I_8=&\gamma\int_{\Gamma_1}|u_t(t)|^{m-2}u_t(t)\int_0^{t}g(t-s)\left(u(t)-u(s)\right){\rm d}s{\rm d}\sigma\nonumber\\
\leq&\frac{\gamma\delta(m-1)}{m}(a-l)\|u_t(t)\|_{m,\Gamma_1}^{m}+\frac{\gamma C_{\ast}^{m}c}{4\delta m}\left(\frac{8pE(0)}{l(p-2)}\right)^{\frac{m-2}{2}}(g\circ\nabla u)(t),
\end {align}
\begin{align}\label{5.19}
I_9=&-\int_{\Gamma_1}u_t(t)\int_0^{t}g'(t-s)\left(u(t)-u(s)\right){\rm d}s{\rm d}\sigma\nonumber\\
\leq&\|u_t(t)\|_{2,\Gamma_1}^2+\frac{a-l}{4\delta\lambda_1}(-g'\circ\nabla u)(t),
\end {align}
\begin{align}\label{5.20}
I_{10}=&-\left(\int_0^{t}g(s){\rm d}s\right)\int_{\Gamma_1}|u_t(t)|^{2}{\rm d}\sigma\leq(a-l)\|u_t(t)\|_{2,\Gamma_1}^2.
\end {align}
A combination of \eqref{5.12}-\eqref{5.20} yields \eqref{5.111}.

{\bf Proof of Theorem \ref{th5.1}.} Since $g$ is continuous and $g(0)>0$, then for any $t_0>0$, we have $$\int_0^{t}g(s){\rm d}s\geq \int_0^{t_0}g(s){\rm d}s:=g_0>0,$$ for all $t\geq t_0$.
From \eqref{4.1}, \eqref{5.6}, \eqref{5.111}, then from \eqref{5.1}, we get
\begin{align}\label{5.21}
L'(t)=&E'(t)+\varepsilon_1 G'(t)+\varepsilon_2 H'(t)\nonumber\\
\leq&-\left[\alpha-\varepsilon_1-\varepsilon_2\left(\alpha+\delta-\int_0^{t}g(s){\rm d}s\right)\right]\|u_t(t)\|_2^2\nonumber\\
&-\left[\frac{1}{2}g(t)-\varepsilon_1\left(-\frac{l}{2}+\frac{C\gamma\delta^{-m}}{m}\left(\frac{2pE(0)}{l(p-2)}\right)^{\frac{m-2}{2}}\right)-\varepsilon_2K_\delta\right]\|\nabla u(t)\|_2^2\nonumber\\
&-\left(\beta-\frac{\beta\varepsilon_2}{2}\right)\|\nabla u_{t}(t)\|_{2}^{2}+\left[\varepsilon_1\frac{(a-l)}{2l}+\varepsilon_2 K_\mu\right]\left(g\circ\nabla u\right)(t)\nonumber\\
&+\left[\varepsilon_1+\varepsilon_2(\delta+a-l)\right]\|u_t(t)\|_{2,\Gamma_1}^2\nonumber\\
&-\left[\gamma-\varepsilon_1\frac{(m-1)\delta^{\frac{m}{m-1}}}{m}-\varepsilon_2\frac{\gamma\delta(m-1)(a-l)}{m}\right]\|u_t(t)\|_{m,\Gamma_1}^{m}\nonumber\\
&+\left[\frac{1}{2}-\varepsilon_2\frac{(a-l)}{2\delta\lambda_1}\right](g'\circ\nabla u)(t)+\varepsilon_1\|u(t)\|_p^p-\frac{\sigma(p-2)}{2p}E(0)E'(t),
\end {align}
where
\begin{align}
K_\delta=\delta+2\delta(a-l)^2+\delta C_{\ast}^{2(p-1)}\left(\frac{2pE(0)}{l(p-2)}\right)^{p-2}>0, \nonumber
\end {align}
\begin{align}
K_\mu=&\frac{a-l}{4\delta}\left(a+\frac{2pbE(0)}{(p-2)l}\right)^2+\frac{\beta(a-l)}{2}+\left(2\delta+\frac{1}{4\delta}\right)(1-l)^2\nonumber\\
&+\frac{\gamma C_{\ast}^{m}c}{4\delta m}\left(\frac{8pE(0)}{l(p-2)}\right)^{\frac{m-2}{2}}+\frac{a-l}{4l}+\frac{a-l}{2\delta\lambda_1}>0.\nonumber
\end {align}
By the trace inequality
\begin{align}
\|u_t\|_{2,\Gamma_1}^2\leq C_1\|u_t\|_{W^{1,2}(\Omega)} \leq C\|u_t\|_2^2,\nonumber
\end {align}
then we have
\begin{align}\label{5.22}
L'(t)\leq&-\left[\alpha-\varepsilon_1-\varepsilon_2\left(\alpha+\delta-\int_0^{t}g(s){\rm d}s\right)-C\left[\varepsilon_1+\varepsilon_2(\delta+a-l)\right]\right]\|u_t(t)\|_2^2\nonumber\\
&-\left[\frac{1}{2}g(t)-\varepsilon_1\left(-\frac{l}{2}+\frac{C\gamma\delta^{-m}}{m}\left(\frac{2pE(0)}{l(p-2)}\right)^{\frac{m-2}{2}}\right)-\varepsilon_2K_\delta\right]\|\nabla u(t)\|_2^2\nonumber\\
&-\left[\beta-\frac{\beta\varepsilon_2}{2}\right]\|\nabla u_{t}(t)\|_{2}^{2}+\left[\varepsilon_1\frac{a-l}{2l}+\varepsilon_2 K_\mu\right]\left(g\circ\nabla u\right)(t)\nonumber\\
&-\left[\gamma-\varepsilon_1\frac{(m-1)\delta^{\frac{m}{m-1}}}{m}-\varepsilon_2\frac{\gamma\delta(m-1)(a-l)}{m}\right]\|u_t(t)\|_{m,\Gamma_1}^{m}\nonumber\\
&+\left[\frac{1}{2}-\varepsilon_2\frac{(a-l)}{2\delta\lambda_1}\right](g'\circ\nabla u)(t)+\varepsilon_1\|u(t)\|_p^p-\frac{\sigma(p-2)}{2p}E(0)E'(t).
\end {align}
At this point, we choose $\delta>0$ satisfying
\begin{align}
0<\delta<\frac{\alpha-\varepsilon_1-\varepsilon_2\left(\alpha-\int_0^{t}g(s){\rm d}s\right)-C\left[\varepsilon_1+\varepsilon_2(a-l)\right]}{(1+C)\varepsilon_2},\nonumber
\end {align}
so we get
\begin{align}\label{5.23}
k_1=\alpha-\varepsilon_1-\varepsilon_2\left(\alpha+\delta-\int_0^{t}g(s){\rm d}s\right)-C\left[\varepsilon_1+\varepsilon_2(\delta+a-l)\right]>0.
\end {align}
We then choose $\varepsilon_1$ and $\varepsilon_2$ so small that \eqref{5.4} and \eqref{5.23} remain valid and
\begin{align}
k_2=\frac{1}{2}g(t)-\varepsilon_1\left(-\frac{l}{2}+\frac{C\gamma\delta^{-m}}{m}\left(\frac{2pE(0)}{l(p-2)}\right)^{\frac{m-2}{2}}\right)-\varepsilon_2K_\delta>0,\nonumber
\end {align}
\begin{align}
k_3=\beta-\frac{\beta\varepsilon_2}{2}>0,\nonumber
\end {align}
\begin{align}
k_4=\gamma-\varepsilon_1\frac{(m-1)\delta^{\frac{m}{m-1}}}{m}-\varepsilon_2\frac{\gamma\delta(m-1)(a-l)}{m}>0,\nonumber
\end {align}
\begin{align}
k_5=\frac{1}{2}-\varepsilon_2\frac{(a-l)}{2\delta\lambda_1}>0.\nonumber
\end {align}
Therefore, we have
\begin{align}\label{5.24}
L'(t)\leq&-k_1\|u_t(t)\|_2^2-k_2\|\nabla u(t)\|_2^2-k_3\|\nabla u_{t}(t)\|_{2}^{2}+\left[\varepsilon_1\frac{a-l}{2l}+\varepsilon_2 K_\mu\right]\left(g\circ\nabla u\right)(t)\nonumber\\
&-k_4\|u_t(t)\|_{m,\Gamma_1}^{m}+k_5(g'\circ\nabla u)(t)+\varepsilon_1\|u(t)\|_p^p-\frac{\sigma(p-2)}{2p}E(0)E'(t).
\end {align}
Then, by $(G1)$ and $(G2)$, for positive constants $M$ and $\alpha_3$ we obtain
\begin{align}\label{5.25}
L'(t)\leq-ME(t)+\alpha_3(g\circ\nabla u)(t)-\frac{\sigma(p-2)}{2p}E(0)E'(t),\ \forall t\geq t_0.
\end {align}
Assume that $(G1)$ and $(G2)$ are satisfied, multiplying \eqref{5.25} by $\xi(t)$, we have
\begin{align}\label{5.26}
\xi(t)L'(t)\leq-M\xi(t)E(t)+\alpha_3\xi(t)(g\circ\nabla u)(t)-\frac{\sigma(p-2)}{2p}E(0)\xi(t)E'(t),\ \forall t\geq t_0,
\end {align}
using $(G2)$, \eqref{4.1} and the fact that $\xi$ and $g$ are nonincreasing, we get
\begin{align}
\xi(t)\int_{0}^{t}g(t-s)\|\nabla u(t)-\nabla u(s)\|_{2}^{2}ds\leq&-\int_{0}^{t}g'(t-s)\|\nabla u(t)-\nabla u(s)\|_{2}^{2}ds\nonumber\\
\leq&-2E'(t)\nonumber
\end {align}
and
\begin{align}
\frac{\sigma(p-2)}{2p}E(0)\xi(t)E'(t)\leq\frac{\sigma(p-2)}{2p}E(0)\xi(0)E'(t).\nonumber\\
\end {align}
Inserting the last two inequalities in \eqref{5.26}, we can obtain
\begin{align}\label{5.27}
\xi(t)L'(t)+2\alpha_4E'(t)\leq-M\xi(t)E(t),\ \forall t\geq t_0,
\end {align}
where $\alpha_4=\alpha_3+\frac{\sigma(p-2)}{2p}E(0)\xi(0)$ is a positive constant.\\
Now we consider the functional
\begin{align}
\Phi=\xi L +2\alpha_4E.\nonumber
\end {align}
Using the facts that $\Phi\sim E$ and $\xi$ is nonincreasing, \eqref{5.27} gives
\begin{align}
\Phi'(t)\leq -k\xi(t)\Phi(t), \ \forall t\geq t_0 ,\nonumber
\end {align}
for some $k>0$. Therefore, direct integration leads to
\begin{align}
\Phi(t)\leq\Phi(t_0)e^{ -k\int_{t_0}^t\xi(s){\rm d}s}, \ \forall t\geq t_0 \nonumber
\end {align}
and the fact that $\Phi\sim E$ yields
\begin{align}
E(t)\leq \alpha_4\Phi(t_0)e^{ -k\int_{t_0}^t\xi(s){\rm d}s}=Ke^{-k\int_{t_0}^t\xi(s){\rm d}s}, \ \forall t\geq t_0, \nonumber
\end {align}
where $\alpha_4$ is a positive constant and $K=\alpha_4\Phi(t_0)$.
This completes the proof.

\section{Blow up} \label{6}
In this section, we prove a finite time blow-up result for initial data in the unstable set. From the definition \eqref{2.4} of the potential depth $d$, for $u\in H_{\Gamma_{0}}^{1}(\Omega)\backslash\{0\}$, we have
\begin{theo}\label{th6.1}  Suppose $2\leq p\leq \overline q$, $m=2$, $(G1)$ and $(G2)$  hold. Let $u$ be the unique local solution to problem \eqref{1.1}, for any fixed $\delta<1$,
assume that $u_0$, $u_1$ satisfy
\begin{align}\label{6.1}
E(0)=\delta d_1,
\end {align}
\begin{align}\label{6.2}
I(0)=0.
\end {align}
Suppose that
\begin{align}\label{6.3}
\int_0^\infty g(s){\rm d}s\leq\frac{p-2}{a(p-2)+1/[(1-\widetilde\delta)^2p+2\delta(1-\widetilde\delta)]},
\end {align}
where $\widetilde\delta =max\{0,\delta\}$ and suppose further that $\int_\Omega u_0 u_1 {\rm d}x+\int_{\Gamma_1} u_0 u_1 {\rm d}\sigma >0$ for $0\leq E(0)<d_1$, then $T<\infty$.
For $t\geq 0$, by \eqref{2.4}, we have $d_1$ is the lower bound of $d(t)$.
\end{theo}
\begin{lemma}\label{Lemma 6.2}
Under the same assumption as in Theorem \ref{th6.1}, one has $I(t)<0$ and
\begin{align}\label{6.4}
d_1<\frac{p-2}{2p}\left[\left(a-\int_0^tg(s)ds\right)\|\nabla u(t)\|_2^2+(g\circ \nabla u)(t)+\frac{b}{2}\|\nabla u(t)\|_2^4\right]<\frac{p-2}{2p}\|u(t)\|_p^p,
\end {align}
for all $t\in [0,T_{\rm max})$.
\end{lemma}

{\bf Proof.} By Lemma \ref{Lemma 4.1} and \eqref{6.1}, we have $E(t)\leq\delta d_1$, for all $t\in[0,T_{\rm max})$. Moreover, we can obtain $I(t)<0$ for all  $t\in[0,T_{\rm max})$. However, if it is not true, then there exists some $t^\ast\in [0,T_{\rm max})$ such that $I(t^\ast)=0$. Thus $I(t)<0$ for all $0\leq t<t^\ast$, i.e.
\begin{align}\label{6.5}
\left(a-\int_0^tg(s)ds\right)\|\nabla u(t)\|_2^2+(g\circ \nabla u)(t)+\frac{b}{2}\|\nabla u(t)\|_2^4<\|u(t)\|_p^p, \quad 0\leq t<t^\ast.
\end {align}
By the definition of $d$, we have
\begin{align}\label{6.6}
d<\frac{p-2}{2p}\left[\left(a-\int_0^tg(s)ds\right)\|\nabla u(t)\|_2^2+(g\circ \nabla u)(t)+\frac{b}{2}\|\nabla u(t)\|_2^4\right], \quad 0\leq t<t^\ast.
\end {align}
Combined with \eqref{6.5} and \eqref{6.6}, we get
\begin{align}
\|u(t)\|_p^p>\frac{2p}{p-2} d_1>0, \quad 0\leq t<t^\ast.\nonumber
\end {align}
By the continuity of $t\mapsto\|u(t)\|_p^p$, we get $u(t^\ast)\neq 0$. Because $d_1$ is the lower bound of $d$ and \eqref{2.2}, we obtain
\begin{align}
d_1\leq\frac{p-2}{2p}\|u(t^\ast)\|_p^p=J(u(t^\ast)),\nonumber
\end {align}
which contradicts to $J(u(t^\ast))\leq E(t^*)<d_1$. By repeating the previous step, we obtain \eqref{6.4}.
\begin{lemma}\label{Lemma 6.3}
Let $F(t)$ be a positive, twice differentiable function, for $t>0$, we can obtain
\begin{align}
F(t)F''(t)-(1+\kappa)F'(t)^{2}\geq 0,\nonumber
\end {align}
with some $\kappa>0$. If $F(0)>0$ and $F'(0)>0$, then there exists a time $T^{\ast}\leq F(0)/[\kappa F'(0)]$ such that $\lim_{t\rightarrow T^{\ast-}}F(t)=\infty$. We can see \cite{LEV1974} for more details.
\end{lemma}

{\bf Proof of Theorem \ref{th6.1}.}
Assume by contradiction that the solution $u$ is global. Then for any $T>0$, let us define the functional $F$ as follows
\begin{align}
F(t)=&\|u(t)\|_{2}^{2}+\|u(t)\|_{2,\Gamma_{1}}^{2}+\frac{\sigma}{2}\int_0^t\|\nabla u(s)\|_2^{4}{\rm d}s+\alpha\int_{0}^{t}\|u(s)\|_{2}^{2}{\rm d}s+\beta\int_{0}^{t}\|\nabla u(s)\|_{2}^{2}{\rm d}s\nonumber\\
&+\gamma\int_{0}^{t}\|u(s)\|_{2,\Gamma_1}^{2}{\rm d}s+(T-t)\left[\alpha\|u_0\|_{2}^{2}+\beta\|\nabla u_0\|_{2}^{2}+\gamma\|u_0\|_{2,\Gamma_1}^{2}+\frac{\sigma}{2}\|\nabla u_0\|_2^4\right]\nonumber\\
&+b(t+T_0)^{2}\nonumber,
\end {align}
where $T$ and $T_0$ are positive constants to be chosen later, $b>0$ if $E(0)<0$ and $b=0$ if $E(0)>0$. Then $F(t)>0$ for all $t\in [0,T]$. Furthermore,
\begin{align}\label{6.7}
F'(t)=&2\int_\Omega u_tu{\rm d}x+2\int_{\Gamma_1}u_tu{\rm d}\sigma+\frac{\sigma}{2}\int_0^t\frac{\rm d}{{\rm d}t}\|\nabla u(s)\|_2^{4}{\rm d}s+2\alpha\int_0^{t}\int_\Omega u_tu{\rm d}x{\rm d}s\nonumber \\
&+2\beta\int_0^{t}\int_\Omega \nabla u_t\cdot\nabla u{\rm d}x{\rm d}s+2\gamma\int_0^{t}\int_{\Gamma_1}u_tu{\rm d}\sigma {\rm d}s+2b(t+T_0),
\end {align}
and
\begin{align}
F''(t)=&2\left[\|u_t(t)\|_2^2-\left(a-\int_{0}^{t}g(s){\rm d}s\right)\|\nabla u(t)\|_{2}^{2}+\|u_t(t)\|_{2,\Gamma_1}^{2}-\frac{b}{2}\|\nabla u(t)\|_2^4\right.\nonumber\\
&\left.-\int_{\Omega}\nabla u(t)\cdot\left(\int_{0}^{t}g(t-s)(\nabla u(t)-\nabla u(s)){\rm d}s\right){\rm d}x+\|u\|_{p}^{p}+b\right].\nonumber
\end {align}
Therefore, using the definition of $F(t)$, we get
\begin{align}
F(t)F''(t)-\frac{p+2}{4}F'(t)^{2}=&F(t)F''(t)+(p+2)\bigg[\eta(t)-\left(F(t)-(T-t)\left(\alpha\|u_{0}(t)\|_{2}^{2}\right.\right.\bigg.\nonumber\\
&\left.\left.+\beta\|\nabla u_{0}(t)\|_{2}^{2}+\gamma\|u_{0}(t)\|_{2,\Gamma_1}^2+\frac{\sigma}{2}\|\nabla u_0\|_2^4\right)\right)\bigg(\|u_t(t)\|_{2}^{2}\bigg.\nonumber\\
&\left.\left.+\|u_t(t)\|_{2,\Gamma_1}^{2}+\frac{\sigma}{2}\int_0^t\|\nabla u(s)\|_2^{4}{\rm d}s+\alpha\int_{0}^{t}\|u_t(s)\|_{2}^{2}{\rm d}s\right.\right.\nonumber\\
&\left.\left.+\beta\int_{0}^{t}\|\nabla u_t(s)\|_{2}^{2}{\rm d}s+\gamma\int_{0}^{t}\|u_t(s)\|_{2,\Gamma_1}^{2}{\rm d}s+b\right)\right]\nonumber
\end {align}
where the function $\eta$ is defined by
\begin{align}
\eta(t)=&\left(\|u(t)\|_{2}^{2}+\|u(t)\|_{2,\Gamma_1}^{2}+\frac{\sigma}{2}\int_0^t\|\nabla u(s)\|_2^{4}{\rm d}s+\alpha\int_{0}^{t}\|u(s)\|_{2}^{2}{\rm d}s+\beta\int_{0}^{t}\|\nabla u(s)\|_{2}^{2}{\rm d}s\right.\nonumber\\
&\left.+\gamma\int_{0}^{t}\|u(s)\|_{2,\Gamma_1}^{2}{\rm d}s+b(t+T_0)^{2}\right)\left(\|u_t(t)\|_{2}^{2}+\|u_t(t)\|_{2,\Gamma_1}^{2}+\frac{\sigma}{2}\int_0^t\|\nabla u_t(s)\|_2^{4}{\rm d}s\right.\nonumber\\
&\left.+\alpha\int_{0}^{t}\|u_t(s)\|_{2}^{2}{\rm d}s+\beta\int_{0}^{t}\|\nabla u_t(s)\|_{2}^{2}{\rm d}s+\gamma\int_{0}^{t}\|u_t(s)\|_{2,\Gamma_1}^{2}{\rm d}s+b\right)-\left[\int_\Omega u_t(t)u(t){\rm d}x\right.\nonumber\\
&\left.+\int_{\Gamma_1}u_t(t)u(t){\rm d}\sigma+\frac{\sigma}{4}\int_0^t\frac{\rm d}{{\rm d}t}\|\nabla u(s)\|_2^{4}{\rm d}s+\alpha\int_0^{t}\int_\Omega u_t(t)u(t){\rm d}x{\rm d}s\right.\nonumber\\
&\left.+\beta\int_0^{t}\int_\Omega \nabla u_t(t)\nabla u(t){\rm d}x{\rm d}s
+\gamma\int_0^{t}\int_{\Gamma_1}u_t(t)u(t){\rm d}\sigma {\rm d}s+b(t+T_0)\right]^{2}\nonumber.
\end {align}
By Cauchy-Schwarz inequality, for all $t\in[0,T]$, we have $\eta(t)\geq0$. As a consequence, we read the following differential inequality
\begin{align}\label{6.8}
F(t)F''(t)-\frac{p+2}{4}F'(t)^{2}\geq F(t)\zeta(t),\quad \forall t\in[t,T].
\end {align}
where
\begin{align}
\zeta(t)=&-p\|u_t(t)\|_{2}^{2}-2\left(a-\int_{0}^{t}g(s){\rm d}s\right)\|\nabla u(t)\|_{2}^{2}-\frac{b}{2}\|\nabla u(t)\|_2^4-p\|u_t(t)\|_{2,\Gamma_1}^{2}+2\|u(t)\|_{p}^{p}\nonumber\\
&-(p+2)\left(\alpha\int_{0}^{t}\|u_t(s)\|_{2}^{2}{\rm d}s+\beta\int_{0}^{t}\|\nabla u_t(s)\|_{2}^{2}{\rm d}s+\gamma\int_{0}^{t}\|u_t(s)\|_{2,\Gamma_1}^{2}{\rm d}s\right.\nonumber\\
&\left.+\frac{\sigma}{2}\int_0^t\|\nabla u_t(s)\|_2^{4}{\rm d}s\right)-2\int_{\Omega}\nabla u(t)\cdot\left(\int_{0}^{t}g(t-s)(\nabla u(t)-\nabla u(s))ds\right){\rm d}x-pb\nonumber\\
=&-2pE(t)+(p-2)\left(a-\int_{0}^{t}g(s){\rm d}s\right)\|\nabla u(t)\|_{2}^{2}+p(g\circ\nabla u)(t)+\frac{(p-1)b}{2}\|\nabla u(t)\|_2^4\nonumber\\
&-(p+2)\left(\alpha\int_{0}^{t}\|u_t(s)\|_{2}^{2}{\rm d}s+\beta\int_{0}^{t}\|\nabla u_t(s)\|_{2}^{2}{\rm d}s+\gamma\int_{0}^{t}\|u_t(s)\|_{2,\Gamma_1}^{2}{\rm d}s\right.\nonumber\\
&\left.+\frac{\sigma}{2}\int_0^t\|\nabla u_t(s)\|_2^{4}{\rm d}s\right)-2\int_{\Omega}\nabla u(t)\cdot\left(\int_{0}^{t}g(t-s)(\nabla u(t)-\nabla u(s)){\rm d}s\right){\rm d}x-pb\nonumber.
\end {align}
From the equality \eqref{4.1}, $(G1)$, $(G2)$ and $m=2$ we have
\begin{align}
E'(t)=&\frac{1}{2}(g'\circ\nabla u)(t)-\frac{1}{2}g(t)\|\nabla u(t)\|_{2}^{2}-\alpha\|u_{t}(t)\|_{2}^{2}-\beta\|\nabla u_{t}(t)\|_{2}^{2}\nonumber\\
&-\gamma\| u_{t}(t)\|_{2,\Gamma_{1}}^{2}-\sigma\left(\frac{1}{2}\frac{\rm d}{\rm dt}\|\nabla u(t)\|_{2}^{2}\right)^2\nonumber\\
\leq&-\alpha\|u_{t}(t)\|_{2}^{2}-\beta\|\nabla u_{t}(t)\|_{2}^{2}-\gamma\| u_{t}(t)\|_{2,\Gamma_{1}}^{2}-\sigma\left(\frac{1}{2}\frac{\rm d}{\rm dt}\|\nabla u(t)\|_{2}^{2}\right)^2\nonumber.
\end {align}
Thus, we can write
\begin{align}\label{6.9}
\zeta(t)\geq&-2pE(0)+(p-2)\left(a-\int_{0}^{t}g(s){\rm d}s\right)\|\nabla u(t)\|_{2}^{2}+p(g\circ\nabla u)(t)+\frac{(p-1)b}{2}\|\nabla u(t)\|_2^4\nonumber\\
&+(p-2)\left(\alpha\int_{0}^{t}\|u_t(s)\|_{2}^{2}{\rm d}s+\beta\int_{0}^{t}\|\nabla u_t(s)\|_{2}^{2}{\rm d}s+\gamma\int_{0}^{t}\|u_t(s)\|_{2,\Gamma_1}^{2}{\rm d}s\right.\nonumber\\
&\left.+\frac{\sigma}{2}\int_0^t\|\nabla u_t(s)\|_2^{4}{\rm d}s\right)-2\int_{\Omega}\nabla u(t)\cdot\left(\int_{0}^{t}g(t-s)(\nabla u(t)-\nabla u(s)){\rm d}s\right){\rm d}x-pb.
\end {align}
By Young's inequality, we have
\begin{align}\label{6.10}
2\int_{\Omega}\nabla u(t)\cdot\left(\int_{0}^{t}g(t-s)(\nabla u(t)-\nabla u(s)){\rm d}s\right)\nonumber\\
\leq\frac{1}{\varepsilon}\int_{0}^{t}g(s)\|\nabla u(t)\|_{2}^{2}{\rm d}s+\varepsilon(g\circ\nabla u)(t),
\end {align}
for any $\varepsilon>0$. Substitute \eqref{6.10} for the fifth term of the right hand side of \eqref{6.9}, we have
\begin{align}\label{6.11}
\zeta(t)\geq&-2pE(0)+\left[a(p-2)-\left(p-2+\frac{1}{\varepsilon}\right)\int_{0}^{t}g(s){\rm d}s\right]\|\nabla u(t)\|_{2}^{2}+(p-\varepsilon)(g\circ\nabla u)(t)\nonumber\\
&+(p-2)\left(\alpha\int_{0}^{t}\|u_t(s)\|_{2}^{2}{\rm d}s+\beta\int_{0}^{t}\|\nabla u_t(s)\|_{2}^{2}{\rm d}s+\gamma\int_{0}^{t}\|u_t(s)\|_{2,\Gamma_1}^{2}{\rm d}s\right.\nonumber\\
&\left.+\frac{\sigma}{2}\int_0^t\|\nabla u_t(s)\|_2^{4}{\rm d}s\right)+\frac{(p-1)b}{2}\|\nabla u(t)\|_2^4-pb.
\end {align}

On the one hand, if $\delta<0$, then $E(0)<0$, we choose $\varepsilon=p$ in \eqref{6.11} and $b$ small enough such that $b\leq-2E(0)$. Then, by \eqref{6.3}, we get
\begin{align}\label{6.12}
\zeta(t)\geq&\left[a(p-2)-\left(p-2+\frac{1}{p}\right)\int_{0}^{t}g(s){\rm d}s\right]\|\nabla u(t)\|_{2}^{2}+p(-2E(0)-b)+\frac{(p-1)b}{2}\|\nabla u(t)\|_2^4\nonumber\\
&+(p-2)\left(\alpha\int_{0}^{t}\|u_t(s)\|_{2}^{2}{\rm d}s+\beta\int_{0}^{t}\|\nabla u_t(s)\|_{2}^{2}{\rm d}s+\gamma\int_{0}^{t}\|u_t(s)\|_{2,\Gamma_1}^{2}{\rm d}s\right.\nonumber\\
&\left.+\frac{\sigma}{2}\int_0^t\|\nabla u_t(s)\|_2^{4}{\rm d}s\right)\nonumber\\
\geq&\left[a(p-2)-\left(p-2+\frac{1}{p}\right)\int_{0}^{t}g(s){\rm d}s\right]\|\nabla u(t)\|_{2}^{2}+\frac{(p-1)b}{2}\|\nabla u(t)\|_2^4\nonumber\\
&+(p-2)\left(\alpha\int_{0}^{t}\|u_t(s)\|_{2}^{2}{\rm d}s+\beta\int_{0}^{t}\|\nabla u_t(s)\|_{2}^{2}{\rm d}s+\gamma\int_{0}^{t}\|u_t(s)\|_{2,\Gamma_1}^{2}{\rm d}s\right.\nonumber\\
&\left.+\frac{\sigma}{2}\int_0^t\|\nabla u_t(s)\|_2^{4}{\rm d}s\right)\geq0.
\end {align}

On the other hand, if $0\leq\delta<1$, then $0\leq E(0)=\delta d_1<d_1$, we choose $\varepsilon=(1-\delta)p+2\delta$ in \eqref{6.11}. Then, we have
\begin{align}
\zeta(t)\geq&-2pE(0)+\delta(p-2)(g\circ\nabla u)(t)+\frac{(p-1)b}{2}\|\nabla u(t)\|_2^4\nonumber\\
&+\left[a(p-2)-\left(p-2+\frac{1}{(1-\delta)p+2\delta}\right)\int_{0}^{t}g(s){\rm d}s\right]\|\nabla u(t)\|_{2}^{2}\nonumber\\
&+(p-2)\left(\alpha\int_{0}^{t}\|u_t(s)\|_{2}^{2}{\rm d}s+\beta\int_{0}^{t}\|\nabla u_t(s)\|_{2}^{2}{\rm d}s+\gamma\int_{0}^{t}\|u_t(s)\|_{2,\Gamma_1}^{2}{\rm d}s\right.\nonumber\\
&\left.+\frac{\sigma}{2}\int_0^t\|\nabla u_t(s)\|_2^{4}{\rm d}s\right).\nonumber
\end {align}
By \eqref{6.3}, we obtain
\begin{align}
a(p-2)-\left(p-2+\frac{1}{(1-\delta)p+2\delta}\right)\int_{0}^{t}g(s){\rm d}s\geq \delta(p-2)\left(a-\int_{0}^{t}g(s){\rm d}s\right),\nonumber
\end {align}
and therefore, by \eqref{6.4}, we have
\begin{align}\label{6.13}
\zeta(t)\geq&-2pE(0)+\delta(p-2)\left[\left(a-\int_{0}^{t}g(s){\rm d}s\right)\|\nabla u(t)\|_{2}^{2}+(g\circ\nabla u)(t)\right]+\frac{(p-1)b}{2}\|\nabla u(t)\|_2^4\nonumber\\
&+(p-2)\left(\alpha\int_{0}^{t}\|u_t(s)\|_{2}^{2}{\rm d}s+\beta\int_{0}^{t}\|\nabla u_t(s)\|_{2}^{2}{\rm d}s+\gamma\int_{0}^{t}\|u_t(s)\|_{2,\Gamma_1}^{2}{\rm d}s\right.\nonumber\\
&\left.+\frac{\sigma}{2}\int_0^t\|\nabla u_t(s)\|_2^{4}{\rm d}s\right)\nonumber\\
\geq&(p-2)\left(\alpha\int_{0}^{t}\|u_t(s)\|_{2}^{2}{\rm d}s+\beta\int_{0}^{t}\|\nabla u_t(s)\|_{2}^{2}{\rm d}s+\gamma\int_{0}^{t}\|u_t(s)\|_{2,\Gamma_1}^{2}{\rm d}s\right.\nonumber\\
&\left.+\frac{\sigma}{2}\int_0^t\|\nabla u_t(s)\|_2^{4}{\rm d}s\right)+\frac{(p-1)b}{2}\|\nabla u(t)\|_2^4+2p(\delta d_1-E(0))\geq0.
\end {align}
From \eqref{6.8}, \eqref{6.12} and \eqref{6.13}, we get
\begin{align}
F(t)F''(t)-\frac{p+2}{4}F'(t)^{2}\geq 0,\quad \forall t\in[t,T].\nonumber
\end {align}
By \eqref{6.7}, We can choose $T_0$ sufficiently large such that $F'(0)=2\int_\Omega u_0u_1{\rm d}x+2\int_{\Gamma_1} u_0u_1{\rm d}\sigma +2bT_0>0$, when $E(0)<0$.
As $0\leq E(0)<d_1$, the condition $\int_\Omega u_0u_1{\rm d}x+\int_{\Gamma_1} u_0u_1{\rm d}\sigma>0$ also ensure that $F'(0)>0$. As $(p+2)/ 4>1$, letting $\kappa=(p-2)/ 4$, by using the concavity argument, we get $\lim_{t\rightarrow T^{\ast-}}F(t)=\infty$, which implies that $\lim_{t\rightarrow T^{*-}}\|\nabla u(t)\|_2^2=\infty$.

%

\end{document}